\documentclass[10pt,onecolumn]{IEEEtran}
\usepackage{amsmath,epsfig}
\usepackage[defs]{ams}
\usepackage[usenames]{color}
\usepackage{umoline}
\usepackage{setspace}
\usepackage{algorithm}
\usepackage{algorithmic}
\usepackage{multirow}
\usepackage{url}

\newcommand{\INPUT}{\item[\textbf{Input:}]}
\newcommand{\OUTPUT}{\item[\textbf{Output:}]}

\newcommand{\norm}[1]{\left\Vert#1\right\Vert}

\newcommand{\supp}[1]{\mathrm{supp}(#1)}

\newcommand{\rank}{\mathrm{rank}}
\newcommand{\atoms}[1]{\mathrm{atoms}(#1)}

\newtheorem{theorem}{Theorem}[section]
\newtheorem{lemma}[theorem]{Lemma}
\newtheorem{proposition}[theorem]{Proposition}
\newtheorem{corollary}[theorem]{Corollary}

\newtheorem{definition}[theorem]{Definition}
\newtheorem{remark}[theorem]{Remark}

\newcommand{\qed}{\nobreak \ifvmode \relax \else
      \ifdim\lastskip<1.5em \hskip-\lastskip
      \hskip1.5em plus0em minus0.5em \fi \nobreak
      \vrule height0.75em width0.5em depth0.25em\fi}

\title{ADMiRA: Atomic Decomposition \\for Minimum Rank Approximation}
\author{Kiryung Lee and Yoram Bresler
\thanks{K. Lee and Y. Bresler are with Coordinated Science Laboratory and Department of ECE, University of Illinois at Urbana-Champaign, IL 61801 USA
e-mail: \{klee81,ybresler\}@illinois.edu}}

\begin{document}
%\doublespacing

\maketitle

\begin{abstract}
We address the inverse problem that arises in compressed sensing of a low-rank matrix.
Our approach is to pose the inverse problem as an approximation problem with a specified target rank of the solution.
A simple search over the target rank then provides the minimum rank solution satisfying a prescribed data approximation bound.
We propose an atomic decomposition that provides an analogy between parsimonious representations of a sparse vector and a low-rank matrix.
Efficient greedy algorithms to solve the inverse problem for the vector case are extended to the matrix case through this atomic decomposition.
In particular, we propose an efficient and guaranteed algorithm named ADMiRA that extends CoSaMP, its analogue for the vector case.
The performance guarantee is given in terms of the rank-restricted isometry property
and bounds both the number of iterations and the error in the approximate solution
for the general case where the solution is approximately low-rank and the measurements are noisy.
With a sparse measurement operator such as the one arising in the matrix completion problem, the computation in ADMiRA is linear in the number of measurements.
The numerical experiments for the matrix completion problem show that,
although the measurement operator in this case does not satisfy the rank-restricted isometry property,
ADMiRA is a competitive algorithm for matrix completion.
\end{abstract}


\section{Introduction}
\label{sec:intro}

Recent studies in compressed sensing have shown that
a sparsity prior in the representation of the unknowns can guarantee unique and stable solutions to underdetermined linear systems.
The idea has been generalized to the matrix case \cite{fazel2008csa} with the rank replacing sparsity to define the parsimony of the representation of the unknowns.
Compressed sensing of a low-rank matrix addresses the inverse problem of reconstructing an unknown low-rank matrix $X_0 \in \mathbb{C}^{m \times n}$
from its linear measurements $b = \mathcal{A} X_0$
\footnote{
Linear operator $\mathcal{A}$ is not a matrix in this equation.
In this paper, to distinguish general linear operators from matrices, we use calligraphic font for general linear operators.
}
via a given linear operator $\mathcal{A}: \mathbb{C}^{m \times n} \rightarrow \mathbb{C}^{p}$.
As in the vector case, the inverse problem is ill-posed in the sense that the number of measurements is much smaller than the number of the unknowns.
Continuing the analogy with  the vector case,
the remarkable fact is that the number of measurements sufficient for unique and stable recovery is roughly on the same order as the number of degrees of freedom in the unknown low rank matrix.
Moreover, under certain conditions, the recovery can be accomplished by polynomial-time algorithms \cite{recht2007gmr}.

One method to solve the inverse problem by exploiting the prior that $X_0$ is low-rank is
to solve the rank minimization problem $\text{P1}$, to minimize the rank within the affine space defined by $b$ and $\mathcal{A}$:
\begin{equation*}
\text{P1:} \qquad
\begin{array}{llll}
\displaystyle \min_{X \in \mathbb{C}^{m \times n}} & \rank(X) \\
\mathrm{subject~to} & \mathcal{A}X = b.
\end{array}
\end{equation*}
In practice, in the presence of measurement noise or modeling error,
a more appropriate measurement model is $b = \mathcal{A} X_0 + \nu$ where the perturbation $\nu$ has bounded Euclidean norm, $\norm{\nu}_2 \leq \eta$.
In this case, the rank minimization problem is written as
\begin{equation*}
\text{P1':} \qquad
\begin{array}{llll}
\displaystyle \min_{X \in \mathbb{C}^{m \times n}} & \rank(X) \\
\mathrm{subject~to} & \norm{\mathcal{A}X - b}_2 \leq \eta
\end{array}
\end{equation*}
with an ellipsoidal constraint.
Indeed, rank minimization has been studied in more general setting where the feasible set is not necessarily restricted as either an affine space or an ellipsoid.
However, due to the non-convexity of the rank, rank minimization is NP-hard even when the feasible set is convex.
Fazel, Hindi, and Boyd \cite{fazel2001rmh} proposed a convex relaxation of the rank minimization problem
by introducing a convex surrogate of $\rank(X)$, which is known as nuclear norm $\norm{X}_*$ and denotes the sum of all singular values of matrix $X$.

Recht, Fazel, and Parrilo \cite{recht2007gmr} studied rank minimization in the framework of compressed sensing
and showed that rank minimization for the matrix case is analogous to $\ell_0$-norm (number of nonzero elements) minimization for the vector case.
They provided an analogy between the two problems and their respective solutions by convex relaxation.
In the analogy, $\ell_1$-norm minimization for the $\ell_0$-norm minimization problem is analogous to nuclear norm minimization for rank minimization.
Both are efficient algorithms, with guaranteed performance under certain conditions, to solve NP-hard problems:
$\ell_0$-norm minimization and rank minimization, respectively.
The respective conditions are given by the sparsity-restricted isometry property \cite{candes2008rip} and the rank-restricted isometry property \cite{recht2007gmr}, \cite{fazel2008csa}, respectively.
However, whereas $\ell_1$-norm minimization corresponds to a linear program (or a quadratically constrained linear program for the noisy case),
nuclear norm minimization is formulated as a convex semidefinite program (SDP).
Although there exist polynomial time algorithms to solve SDP, in practice they do not scale well to large problems.

Recently, several authors proposed methods for solving large scale SDP derived from rank minimization.
These include interior point methods for SDP, projected subgradient methods, and low-rank parametrization \cite{recht2007gmr} combined with a customized interior point method \cite{liu2008ipm}.
These methods can solve larger rank minimization problems, which the general purpose SDP solvers cannot.
However, the dimension of the problem is still restricted and some of these methods do not guarantee convergence to the global minimum.
Cai, Candes, and Shen \cite{cai2008svt} proposed singular value thresholding (SVT), which penalizes the objective of nuclear norm minimization by the squared Frobenius norm.
The dual of the penalized problem admits a projected subgradient method where the updates can be done by computing truncated singular value decompositions.
They have shown that the solution given by SVT converges to the solution to nuclear norm minimization as the penalty parameter increases.
However, an analysis of the convergence rate is missing and hence the quality of the solution obtained by this method is not guaranteed.
Furthermore, the efficiency of SVT is restricted to the noiseless case where the constraint is affine (\textit{i.e.}, linear equality).
Ma, Goldfarb, and Chen \cite{ma2008fpb} proposed a formulation of nuclear norm minimization by using the Bregman divergence
that admits an efficient fixed point algorithm, which is also based on the singular value decomposition.
They did not provide a convergence rate analysis and the efficiency is also restricted to the noiseless, affine constraint case.
Meka \textit{et. al.} \cite{meka2008rmv} used multiplicative updates and online convex programming to provide an approximate solution to rank minimization.
However, their result depends on the (unverified) existence of an oracle that provides the solution to the rank minimization problem with a single linear constraint in constant time.

An alternative method to solve the inverse problem of compressed sensing of a matrix is minimum rank approximation,
\begin{equation*}
\text{P2:} \qquad
\begin{array}{llll}
\displaystyle \min_{X \in \mathbb{C}^{m \times n}} & \norm{\mathcal{A}X - b}_2 \\
\mathrm{subject~to} & \rank(X) \leq r,
\end{array}
\end{equation*}
where $r = \rank(X_0)$ denotes the minimum rank.
The advantage of formulation $\text{P2}$ is that it can handle both the noiseless case and the noisy case in a single form.
It also works for more general case where $X_0$ is not exactly low-rank but admits an accurate approximation by a low-rank matrix.
When the minimum rank $r$ is unknown, an incremental search over $r$ will increase the complexity of the solution by at most factor $r$.
If an upper bound on $r$ is available, then a bisection search over $r$ can be used because the minimum of $\text{P2}$ is monotone decreasing in $r$.
Hence the factor reduces to $\log r$.
Indeed, this is not an issue in many applications where the rank is assumed to be a small constant.

Recently, several algorithms have been proposed to solve $\text{P2}$.
Halder and Diego \cite{haldar2009rcs} proposed an alternating least square approach by exploiting the explicit factorization of a rank-$r$ matrix.
Their algorithm is computationally efficient but does not provide any performance guarantee.
Keshavan, Oh, and Montanari \cite{keshavan2009mcf} proposed an algorithm based on optimization over the Grassmann manifold.
Their algorithm first finds a good starting point by an operation called trimming and minimizes the objective of $\text{P2}$ using a line search and gradient descent over the Grassmann manifold.
They provide a performance guarantee only for the matrix completion problem where the linear operator $\mathcal{A}$ takes a few entries from $X_0$.
Moreover, the performance guarantee is restricted to the noiseless case.

Minimum rank approximation, or rank-$r$ approximation for the matrix case, is analogous to $s$-term approximation for the vector case.
Like rank-$r$ matrix approximation, $s$-term vector approximation is a way to find the sparsest solution of an ill-posed inverse problem in compressed sensing.
For $s$-term approximation, besides efficient greedy heuristics such as Matching Pursuit (MP) \cite{mallat1993mpt} and Orthogonal Matching Pursuit (OMP) \cite{pati1993omp},
there are recent algorithms, which are more efficient than convex relaxation and also have performance guarantees.
These include Compressive Sampling Matching Pursuit (CoSaMP) \cite{needell2009cis} and Subspace Pursuit (SP) \cite{dai2008spc}.
To date, no such algorithms have been available for the matrix case.

In this paper, we propose an iterative algorithm for the rank minimization problem, which is a generalization
\footnote{
There is another generalization of CoSaMP, namely model-based CoSaMP \cite{baraniuk2008mbc}.
However, this generalization addresses a completely different and unrelated problem: sparse vector approximation subject to a special (e.g., tree) structure.
Furthermore, the extensions of CoSaMP to model-based CoSaMP and to ADMiRA are independent:
neither one follows from the other, and neither one is a special case of the other.
}
of the CoSaMP algorithm to the matrix case.
We call this algorithm ``Atomic Decomposition for Minimum Rank Approximation,'' abbreviated as ADMiRA.
ADMiRA is computationally efficient in the sense that
the core computation consists of least squares and truncated singular value decompositions,
which are both basic linear algebra problems and admit efficient algorithms.
Indeed, ADMiRA is the first guaranteed algorithms among those proposed to solve minimum rank approximation
\footnote{
ADMiRA \cite{leebre2009egr} was followed by the algorithm by Keshavan et. al. \cite{keshavan2009mcf}.
This short version \cite{leebre2009egr} will be presented at ISIT'09.
}.
Furthermore, ADMiRA provides a strong performance guarantee for $\text{P2}$ that covers the general case where $X_0$ is only approximately low-rank and $b$ contains noise.
The strong performance guarantee of ADMiRA is comparable to that of nuclear norm minimization in \cite{fazel2008csa}.
In the noiseless case, SVT \cite{cai2008svt} may be considered a competitor to ADMiRA.
However, for the noisy case, SVT involves more than the simple singular value thresholding operation.

Matrix completion is a special case of low-rank matrix approximation from linear measurements
where the linear operator takes a few random entries of the unknown matrix.
It has received considerable attention owing to its important applications such as collaborative filtering.
However, the linear operator in matrix completion does not satisfy the rank-restricted isometry property \cite{candes2008emc}.
Therefore, at the present time, ADMiRA does not have a guarantee for matrix completion.
None the less, empirical performance on matrix completion is better than SVT (for the experiments in this paper).

The remaining of this paper is organized as follows:
The atomic decomposition and the analogy between the greedy algorithm for the vector case and the matrix case are introduced in Section~\ref{sec:vector_vs_matrix}.
The new algorithm ADMiRA and its performance guarantee are explained in Section~\ref{sec:algorithm} and in Section~\ref{sec:performance_guarantee}, respectively.
By using the tools in Section~\ref{sec:properties_rrip}, the performance guarantees are derived in Section~\ref{sec:proof_performance_guarantee} and Section~\ref{sec:iter_cnt}.
Implementation issues and the computational complexity are discussed in Section~\ref{sec:large}
and numerical results in Section~\ref{sec:num_exp}, followed by conclusions.
Our exposition of ADMiRA follows the line of Needell and Tropp's exposition of CoSaMP \cite{needell2009cis},
to highlight, on the one hand, the close analogy, and on the other hand the differences between the two algorithms and their analysis.
Indeed, there exist significant differences between rank-$r$ approximation for the matrix case and $s$-term approximation for the vector case, which are discussed in some detail.

\section{Vector vs Matrix}
\label{sec:vector_vs_matrix}

\subsection{Preliminaries}

Throughout this paper, we use two vector spaces:
the space of column vectors $\mathbb{C}^p$ and the space of matrices $\mathbb{C}^{m \times n}$.
For $\mathbb{C}^p$, the inner product is defined by
$\langle x , y \rangle_{\mathbb{C}^p} = y^H x$ for $x, y \in \mathbb{C}^p$ where $y^H$ denotes the Hermitian transpose of $y$,
and the induced Hilbert-Schmidt norm is the Euclidean or $\ell_2$-norm given by
$\norm{x}_2^2 = \langle x, x \rangle_{\mathbb{C}^p}$ for $x \in \mathbb{C}^p$.
For $\mathbb{C}^{m \times n}$, the inner product is defined by
$\langle X, Y \rangle_{\mathbb{C}^{m \times n}} = \mathrm{tr}(Y^H X)$ for $X, Y \in \mathbb{C}^{m \times n}$,
and the induced norm is the Frobenius norm given by
$\norm{X}_F^2 = \langle X, X \rangle_{\mathbb{C}^{m \times n}}$ for $X \in \mathbb{C}^{m \times n}$.

\subsection{Atomic Decomposition}

Let $\Gamma$ denote the set of all nonzero rank-one matrices in $\mathbb{C}^{m \times n}$.
We can refine $\Gamma$ so that any two distinct elements are not collinear.
The resulting subset $\mathbb{O}$ is referred to as the \textit{set of atoms}
\footnote{The ``atom'' in this paper is different from Mallat and Zhang's ``atom'' \cite{mallat1993mpt},
which is an element in the dictionary, a finite set of vectors.
In our terminology, an atom is a rank-one matrix, an element in an infinite set of vectors (in the vector space $\mathbb{C}^{m \times n}$).
In both cases, however, an atom denotes an irreducible quantity -- a singleton subset, not representable with fewer elements.
(Indeed, for each atom $\psi$, the corresponding atomic space $\mathrm{span}(\psi)$ is irreducible.)
}
of $\mathbb{C}^{m \times n}$.
Then the \textit{set of atomic spaces} $\mathbb{A}$ of $\mathbb{C}^{m \times n}$ is defined by
$\mathbb{A} \triangleq \{\mathrm{span}(\psi) :~ \psi \in \mathbb{O}\}$.
Each subspace $S \in \mathbb{A}$ is one-dimensional and hence is irreducible in the sense that
$S = S_1 + S_2$ for some $S_1,S_2 \in \mathbb{A}$ implies $S_1 = S_2 = S$.
Since $\mathbb{O}$ is an uncountably infinite set in a finite dimensional space $\mathbb{C}^{m \times n}$,
the elements in $\mathbb{O}$ are not linearly independent.
Regardless of the choice of $\mathbb{O}$, $\mathbb{A}$ is uniquely determined.
Without loss of generality, we fix $\mathbb{O}$ such that all elements have unit Frobenius norm.

Given a matrix $X \in \mathbb{C}^{m \times n}$,
its representation $X = \sum_j \alpha_j \psi_j$ as a linear combination of atoms
is referred to as an \textit{atomic decomposition} of $X$.
Since $\mathbb{O}$ spans $\mathbb{C}^{m \times n}$, an atomic decomposition of $X$ exists for all $X \in \mathbb{C}^{m \times n}$.
A subset $\Psi = \{ \psi \in \mathbb{O} :~ \langle \psi_j, \psi_k \rangle_{\mathbb{C}^{m \times n}} = \delta_{jk} \}$
of unit-norm and pairwise orthogonal atoms in $\mathbb{O}$ will be called an \textit{orthonormal set of atoms}.

\begin{definition}
Let $\mathbb{O}$ be a set of atoms of $\mathbb{C}^{m \times n}$.
Given $X \in \mathbb{C}^{m \times n}$, we define $\atoms{X}$ as the smallest set of atoms in $\mathbb{O}$ that spans $X$,
\begin{equation}
\atoms{X} \triangleq \displaystyle \arg\min_\Psi \left\{ |\Psi| :~ \Psi \subset \mathbb{O}, \quad X \in \mathrm{span}(\Psi) \right\}.
\label{eq:def_atom}
\end{equation}
\end{definition}
Note that $\atoms{X}$ is not unique.

An orthonormal set $\atoms{X} \subset \mathbb{O}$ is given by the singular value decomposition of $X$.
Let $X = \sum_{k=1}^{\rank(X)} \sigma_k u_k v_k^H$ denote the singular value decomposition of $X$ with singular values in decreasing order.
While $u_k v_k^H$ need not be in $\mathbb{O}$, for each $k$, there exists $\rho_k \in \mathbb{C}$ such that $|\rho_k| = 1$ and $\rho_k u_k v_k^H \in \mathbb{O}$.
Then an orthonormal set $\atoms{X} \subset \mathbb{O}$ is given by
\begin{equation*}
\atoms{X} = \{ \rho_k u_k v_k^H \}_{k=1}^{\rank(X)}.
\end{equation*}

\begin{remark}
$\atoms{X}$ and $\rank(X) = |\atoms{X}|$ of a matrix $X \in \mathbb{C}^{m \times n}$ are the counterparts of
$\supp{x}$ and $\norm{x}_0 = |\supp{x}|$ for a vector $x \in \mathbb{C}^p$, respectively.
\end{remark}

\subsection{Generalized Correlation Maximization}

Recht, Fazel, and Parrilo \cite{recht2007gmr} showed an analogy between rank minimization $\text{P1}$ and $\ell_0$-norm minimization.
We consider instead the rank-$r$ matrix approximation problem $\text{P2}$ and its analogue -- the $s$-term vector approximation problem
\begin{equation*}
\text{P3:} \qquad
\begin{array}{llll}
\displaystyle \min_{x \in \mathbb{C}^{n}} & \norm{Ax - b}_2 \\
\mathrm{subject~to} & \norm{x}_0 \leq s.
\end{array}
\end{equation*}
In Problem $\text{P3}$, variable $x$ lives in the union of $s$ dimensional subspaces of $\mathbb{C}^{n}$,
each spanned by $s$ elements in the finite set $\mathbb{E} = \{ e_1, \ldots, e_n \}$, the standard basis of $\mathbb{C}^n$.
Thus the union contains all $s$-sparse vectors in $\mathbb{C}^n$.
Importantly, finitely many ($n \choose s$, to be precise) subspaces participate in the union.
Therefore, it is not surprising that $\text{P3}$ can be solved exactly by exhaustive enumeration,
and finite selection algorithms such as CoSaMP are applicable.

In the rank-$r$ matrix approximation problem $\text{P2}$,
the matrix variable $X$ lives in the union of subspaces of $\mathbb{C}^{m \times n}$,
each of which is spanned by $r$ atoms in the set $\mathbb{O}$.
Indeed, if $X \in \mathbb{C}^{m \times n}$ is spanned by $r$ atoms in $\mathbb{O}$, then
$\rank(X) \leq r$ by the subadditivity of the rank.
Conversely, if $\rank(X) = r$,
then $X$ is a linear combination of rank-one matrices and hence there exist $r$ atoms that span $X$.
Note that uncountably infinitely many subspaces participate in the union.
Therefore, some selection rules in the greedy algorithms for $\ell_0$-norm minimization and $s$-term vector approximation do not generalize in a straightforward way.
None the less, using our formulation of the rank-$r$ matrix approximation problem in terms of an atomic decomposition,
we extend the analogy between the vector and matrix cases,
and propose a way to generalize these selection rules to the rank-$r$ matrix approximation problem.

First, consider the correlation maximization in greedy algorithms for the vector case.
Matching Pursuit (MP) \cite{mallat1993mpt} and Orthogonal Matching Pursuit (OMP) \cite{pati1993omp} choose the index $k \in \{1,\ldots,n\}$
that maximizes the correlation $\left|a_k^H (b - A \hat{x})\right|$ between the $k$-th column $a_k$ of $A$ and the residual in each iteration,
where $\hat{x}$ is the solution of the previous iteration.
Given a set $\Psi$, let $\mathcal{P}_\Psi$ denote the (orthogonal) projection operator onto the subspace spanned by $\Psi$ in the corresponding embedding space.
When $\Psi = \{ \psi \}$ is a singleton set, $\mathcal{P}_\psi$ will denote $\mathcal{P}_\Psi$.
For example, $\mathcal{P}_{e_k}$ denotes the projection operator onto the subspace in $\mathbb{C}^n$ spanned by $e_k$.
From
\begin{eqnarray*}
\left| a_k^H (b - A \hat{x}) \right|
= \left| \langle A^H(b - A \hat{x}), e_k \rangle_{\mathbb{C}^n} \right|
= \norm{\mathcal{P}_{e_k} A^H (b - A \hat{x})}_2,
\end{eqnarray*}
it follows that
maximizing the correlation implies
maximizing the norm of the projection of the image under $A^H$ of the residual $b - A \hat{x}$ onto the selected one dimensional subspace.

The following selection rule generalizes the correlation maximization to the matrix case.
We maximize the norm of the projection over all one-dimensional subspaces spanned by an atom in $\mathbb{O}$:
\begin{equation}
\max_{\psi \in \mathbb{O}} \left| \langle b - \mathcal{A}\widehat{X}, \mathcal{A}\psi \rangle_{\mathbb{C}^{m \times n}} \right|
= \max_{\psi \in \mathbb{O}} \norm{\mathcal{P}_\psi \mathcal{A}^*(b - \mathcal{A}\widehat{X})}_F,
\label{eq:sel_one_largest_cor}
\end{equation}
where $\mathcal{A}^* : \mathbb{C}^p \rightarrow \mathbb{C}^{m \times n}$ denotes the adjoint operator of $\mathcal{A}$.
By the Eckart-Young Theorem, the basis of the best subspace is obtained from the singular value decomposition of $M = \mathcal{A}^*(b - \mathcal{A}\widehat{X})$,
as $\psi = u_1 v_1^H$, where $u_1$ and $v_1$ are the principal left and right singular vectors.
\begin{remark}
Applying the selection rule (\ref{eq:sel_one_largest_cor}) to update $\widehat{X}$ recursively leads to greedy algorithms generalizing MP and OMP to rank minimization.
\end{remark}

Next, consider the rule in recent algorithms such as CoSaMP and SP.
The selection rule chooses the subset $J$ of $\{1,\ldots,n\}$ with $|J| = s$ defined by
\begin{equation}
\left|a_k^H (b - A \hat{x})\right| \geq \left|a_j^H (b - A \hat{x})\right|, \quad \forall k \in J, \forall j \not\in J.
\label{eq:sel_largest_cor}
\end{equation}
This is equivalent to maximizing
\begin{eqnarray*}
\sum_{k \in J} \left|a_k^H (b - A \hat{x})\right|^2
= \sum_{k \in J} \norm{\mathcal{P}_{e_k} A^H (b - A \hat{x})}_2^2
= \norm{\mathcal{P}_{\{e_k\}_{k \in J}} A^H (b - A \hat{x})}_2^2.
\end{eqnarray*}
In other words, selection rule (\ref{eq:sel_largest_cor}) finds the best subspace spanned by $s$ elements in $\mathbb{E}$
that maximizes the norm of the projection of $M = A^H (b - A \hat{x})$ onto that $s$-dimensional subspace.

The following selection rule generalizes the selection rule (\ref{eq:sel_largest_cor}) to the matrix case.
We maximize the norm of the projection over all subspaces spanned by a subset with at most $r$ atoms in $\mathbb{O}$:
\begin{equation*}
\max_{\Psi \subset \mathbb{O}} \left\{ \norm{\mathcal{P}_\Psi \mathcal{A}^*(b - \mathcal{A}\widehat{X})}_F :~ |\Psi| \leq r \right\}
\end{equation*}
A basis $\Psi$ of the best subspace is again obtained from the singular value decomposition of $M = \mathcal{A}^*(b - \mathcal{A}\widehat{X})$,
as $\Psi = \{\rho_k u_k v_k^H\}_{k=1}^r$, where $u_k$ and $v_k$, $k = 1,\ldots,r$ are the $r$ principal left and right singular vectors, respectively
and for each $k$, $\rho_k \in \mathbb{C}$ satisfies $|\rho_k| = 1$
\footnote{
Once the best subspace is determined, it is not required to compute the constants $\rho_k$'s.
}
.
Note that $\Psi$ is an orthonormal set although this is not enforced as an explicit constraint in the maximization.

\section{Algorithm}
\label{sec:algorithm}

\begin{algorithm}
\caption{ADMiRA}
\begin{algorithmic}[1]
\INPUT $\mathcal{A}: \mathbb{C}^{m \times n} \rightarrow \mathbb{C}^p$, $b \in \mathbb{C}^p$, and target rank $r \in \mathbb{N}$
\OUTPUT rank-$r$ solution $\widehat{X}$ to $\text{P2}$
\STATE $\widehat{X} \leftarrow 0$
\STATE $\widehat{\Psi} \leftarrow \emptyset$
\WHILE{stop criterion is false}
\STATE $\Psi' \leftarrow \displaystyle \arg\max_{\Psi \subset \mathbb{O}} \left\{ \norm{\mathcal{P}_\Psi \mathcal{A}^* (b - \mathcal{A}\widehat{X})}_F :~ |\Psi| \leq 2r \right\}$
\label{step:cor_max_A}
\STATE $\widetilde{\Psi} \leftarrow \Psi' \cup \widehat{\Psi}$
\label{step:merge}
\STATE $\widetilde{X} \leftarrow \displaystyle \arg\min_X \left\{ \norm{b - \mathcal{A}X}_2 :~ X \in \mathrm{span}(\widetilde{\Psi}) \right\}$
\label{step:solve_LS}
\STATE $\widehat{\Psi} \leftarrow \displaystyle \arg\max_{\Psi \subset \mathbb{O}} \left\{ \norm{\mathcal{P}_\Psi \widetilde{X}}_F :~ |\Psi| \leq r \right\}$
\label{step:cor_max_C}
\STATE $\widehat{X} \leftarrow \mathcal{P}_{\widehat{\Psi}} \widetilde{X}$
\label{step:proj}
\ENDWHILE
\RETURN $\widehat{X}$
\end{algorithmic}
\label{alg:ADMiRA}
\end{algorithm}

Algorithm~\ref{alg:ADMiRA} describes ADMiRA.
Intuitively, ADMiRA iteratively refines the pair $(\widehat{\Psi},\widehat{X}) \in \mathbb{O} \times \mathbb{C}^{m \times n}$
where $\widehat{\Psi}$ is the set of $r$ atoms that spans an approximate solution $\widehat{X}$ to $\text{P2}$.
Step~\ref{step:cor_max_A} finds a set of $2r$ atoms $\Psi'$ that spans a good approximation of $X_0 - \widehat{X}$,
which corresponds to the information not explained by the solution $\widehat{X}$ in the previous iteration.
Here ADMiRA assumes that $\mathcal{A}$ acts like an isometry on a low-rank matrix $X_0 - \widehat{X}$,
which implies that $\mathcal{A}^* \mathcal{A}$ acts like a (scaled) identity operator on $X_0 - \widehat{X}$.
Under this assumption, the $2r$ leading principal components of the proxy matrix $\mathcal{A}^* (b - \mathcal{A}\widehat{X}) = \mathcal{A}^* \mathcal{A} (X_0 - \widehat{X})$ are a good choice for $\Psi'$.
The quality of a linear approximation of $X_0 - \widehat{X}$ spanned by $\Psi'$ improves as iteration goes.
This will be quantitatively analyzed in the proof of the performance guarantee.
If $\widehat{\Psi}$ and $\Psi'$ span good approximations of $\widehat{X}$ and $X_0 - \widehat{X}$, respectively,
then $\widetilde{\Psi} = \widehat{\Psi} \cup \Psi'$ will span a good approximation of $X_0$.
Steps~\ref{step:solve_LS} and \ref{step:cor_max_C} refine the set $\widetilde{\Psi}$ into a set of $r$ atoms.
We first compute a rank-$3r$ approximate solution $\widetilde{X}$ and then take its best rank-$r$ approximation to get a feasible solution $\widehat{X}$ with rank $r$.
In the process, the set $\widetilde{\Psi}$ of $3r$ atoms is also trimmed to the $r$ atom set $\widehat{\Psi}$ so that it can span an approximate solution $\widehat{X}$ closer to $X_0$.

ADMiRA is guaranteed to converge to the global optimum in at most $6(r+1)$ iterations when the assumptions of ADMiRA in Section~\ref{sec:performance_guarantee} are satisfied.
However, similarly to the vector case \cite{needell2009cis},
it is more difficult to verify the satisfiability of the assumptions than solve the recovery problem itself,
and to date there is no known algorithm to perform this verification.
Instead of relying on the theoretical bound on the number of iterations, we use an empirical stopping criterion below.
If either the monotone decrease of $\|b - \mathcal{A} \widehat{X}\|_2 / \norm{b}_2$ is broken or $\|b - \mathcal{A} \widehat{X}\|_2 / \norm{b}_2$ falls a given threshold, ADMiRA stops.

In terms of computation, Steps~\ref{step:cor_max_A} and \ref{step:cor_max_C} involve finding a best rank-$2r$ or rank-$r$ approximation to a given matrix (e.g., by truncating the SVD),
while Step~\ref{step:solve_LS} involves the solution of a linear least-squares problem -- all standard numerical linear algebra problems.
Step~\ref{step:merge} merges two given sets of atoms in $\mathbb{O}$ by taking their union.
As described in more detail in Section~\ref{sec:large}, these computations can be further simplified and their cost reduced
by storing and operating on the low rank matrices in factored form, and taking advantage of special structure of the measurement operator $\mathcal{A}$, such as sparsity.

Most steps of ADMiRA are similar to those of CoSaMP except Step~\ref{step:cor_max_A} and Step~\ref{step:cor_max_C}.
The common feasible set $\mathbb{O}$ of the maximization problems in Step~\ref{step:cor_max_A} and Step~\ref{step:cor_max_C} is infinite and not orthogonal,
whereas the analogous set $\mathbb{E}$ in CoSaMP is finite and orthonormal.
As a result, the maximization problems over the infinite set $\mathbb{O}$ in ADMiRA are more difficult than those in the analogous steps of CoSaMP,
which can be simply solved by selecting the coordinates with the largest magnitudes.
None the less, singular value decomposition can solve the maximization problems over the infinite set efficiently.

\section{Main Results: Performance Guarantee}
\label{sec:performance_guarantee}

\subsection{Rank-Restricted Isometry Property (R-RIP)}

Recht \textit{et al} \cite{recht2007gmr} generalized the sparsity-restricted isometry property (RIP) defined for sparse vectors to low rank matrices.
They also demonstrated ``nearly isometric families'' satisfying this R-RIP (with overwhelming probability).
These include random linear operators generated from i.i.d. Gaussian, or i.i.d. symmetric Bernoulli distributions.
In order to draw the analogy with known results in $\ell_0$-norm minimization, we slightly modify their definition by squaring the norm in the inequality.
Given a linear operator $\mathcal{A}: \mathbb{C}^{m \times n} \rightarrow \mathbb{C}^p$,
the rank-restricted isometry constant $\delta_r(\mathcal{A})$ is defined as the minimum constant that satisfies
\begin{equation}
(1 - \delta_r(\mathcal{A})) \norm{X}_F^2 \leq \norm{\gamma \mathcal{A} X}_2^2 \leq (1 + \delta_r(\mathcal{A})) \norm{X}_F^2,
\label{eq:rip}
\end{equation}
for all $X \in \mathbb{C}^{m \times n}$ with $\rank(X) \leq r$ for some constant $\gamma > 0$.
Throughout this paper, we assume that the linear operator $\mathcal{A}$ is scaled appropriately so that $\gamma = 1$ in (\ref{eq:rip})
\footnote{
If $\gamma \neq 1$, then the noise term in (\ref{eq:unrecoverable_energy}) needs to be scaled accordingly.
}
.
If $\mathcal{A}$ has a small rank-restricted isometry constant $\delta_r(\mathcal{A}) \ll 1$,
then (\ref{eq:rip}) implies that $\mathcal{A}$ acts like an isometry (scaled by $\gamma$) on the matrices whose rank is equal to or less than $r$.
In this case, $\mathcal{A}$ is called a rank-restricted isometry to indicate that
the domain where $\mathcal{A}$ is nearly an isometry is restricted to the set of low-rank matrices.

\subsection{Performance Guarantee}

Subject to the R-RIP, the Atomic Decomposition for Minimum Rank Approximation Algorithm (ADMiRA) has a performance guarantee analogous to that of CoSaMP.

The followings are the assumptions in ADMiRA:
\begin{description}
\item[A1:] The target rank is fixed as $r$.
\item[A2:] The linear operator $\mathcal{A}$ satisfies $\delta_{4r}(\mathcal{A}) \leq 0.04$.
\item[A3:] The measurement is obtained by
\begin{equation}
b = \mathcal{A} X_0 + \nu,
\label{eq:noisy_meas_mdl}
\end{equation}
where $\nu$ is the discrepancy between the measurement and the linear model $\mathcal{A}X_0$.
No assumptions are made about the matrix $X_0$ underlying the measurement, and it can be arbitrary.
\end{description}

Assumption A2 plays a key role in deriving the performance guarantee of ADMiRA:
it enforces the rank-restricted isometry property of the linear operator $\mathcal{A}$.
Although the verification of the satisfiability of A2 is as difficult as or more difficult than the recovery problem itself,
as mentioned above, nearly isometric families that satisfy the condition in A2 have been demonstrated \cite{recht2007gmr}.

The performance guarantees are specified in terms of a measure of inherent approximation error, termed \emph{unrecoverable energy} defined by
\begin{equation}
\epsilon = \norm{X_0 - X_{0,r}}_F + \frac{1}{\sqrt{r}} \norm{X_0 - X_{0,r}}_* + \norm{\nu}_2,
\label{eq:unrecoverable_energy}
\end{equation}
where $X_{0,r}$ denotes the best rank-$r$ approximation of $X_0$.
The first two terms in $\epsilon$ define a metric of the minimum distance
between the ``true'' matrix $X_0$ and a rank-$r$ matrix.
This is analogous to the notion of a measure of compressibility of a vector in sparse vector approximation.
By the Eckart-Young-Mirsky theorem \cite{mirsky1960sgf}, no rank-$r$ matrix can come closer to $X_0$ in this metric.
In particular, the optimal solution to $\text{P2}$ cannot come closer to $X_0$ in this metric.
The third term is the norm of the measurement noise, which must also limit the accuracy of the approximation provided by a solution to $\text{P2}$.

\begin{theorem}
Let $\widehat{X}_k$ denote the estimate of $X_0$ in the $k$-th iteration of ADMiRA.
For each $k \geq 0$, $\widehat{X}_k$ satisfies the following recursion:
\begin{equation*}
\|X_0 - \widehat{X}_{k+1}\|_F \leq 0.5 \|X_0 - \widehat{X}_k\|_F + 8 \epsilon,
\end{equation*}
where $\epsilon$ is the unrecoverable energy.
From the above relation, it follows that
\begin{equation*}
\|X_0 - \widehat{X}_k\|_F \leq 2^{-k} \norm{X_0}_F + 16 \epsilon, \quad \forall k \geq 0.
\end{equation*}
\label{thm:pg_gen}
\end{theorem}

Theorem~\ref{thm:pg_gen} shows the geometric convergence of ADMiRA.
In fact, convergence in a finite number of steps can be achieved as stated by the following theorem.

\begin{theorem}
After at most $6(r+1)$ iterations, ADMiRA provides a rank-$r$ approximation $\widehat{X}$ of $X_0$, which satisfies
\begin{equation*}
\|X_0 - \widehat{X}\|_F \leq 17 \epsilon,
\end{equation*}
where $\epsilon$ is the unrecoverable energy.
\label{thm:iter_cnt_gen}
\end{theorem}

Depending on the spectral properties of the matrix $X_0$, even faster convergence is possible (See Section~\ref{sec:iter_cnt} for details).

\subsection{Relationship between $\text{P1}$, $\text{P2}$, and ADMiRA}

The approximation $\widehat{X}$ given by ADMiRA is a solution to $\text{P2}$.
When there is no noise in the measurement, \textit{i.e.}, $b = \mathcal{A} X_0$, where $X_0$ is the solution to $\text{P1}$,
Theorem~\ref{thm:pg_gen} states that if the ADMiRA assumptions are satisfied with $r \geq \rank(X_0)$, then $\widehat{X} = X_0$.
An appropriate value can be assigned to $r$ by an incremental search over $r$.

For the noisy measurement case, the linear constraint in $\text{P1}$ is replaced by a quadratic constraint and
the rank minimization problem is written as:
\begin{equation*}
\text{P1$'$:} \qquad
\begin{array}{llll}
\displaystyle \min_{X \in \mathbb{C}^{m \times n}} & \rank(X) \\
\mathrm{subject~to} & \norm{\mathcal{A}X - b}_2 \leq \eta.
\end{array}
\end{equation*}
Let $X'$ denote a minimizer to $\text{P1$'$}$.
In this case, the approximation $\widehat{X}$ produced by ADMiRA is not necessarily equivalent to $X'$,
but by Theorem~\ref{thm:pg_gen} the distance between the two is bounded by $\|X' - \widehat{X}\|_F \leq 17 \eta$
for all $r \geq \rank(X')$ that satisfies the ADMiRA assumptions.

\section{Properties of the Rank-Restricted Isometry}
\label{sec:properties_rrip}

We introduce and prove a number of properties of the rank-restricted isometry.
These properties serve as key tools for proving the performance guarantees for ADMiRA in this paper.
These properties further extend the analogy between the sparse vector and the low-rank matrix approximation problems
($\text{P3}$ and $\text{P2}$, respectively), and are therefore also of interest in their own right.
The proofs are contained in the Appendix.

\begin{proposition}
The rank-restricted isometry constant $\delta_r(\mathcal{A})$ is nondecreasing in $r$.
\label{prop:monotone_ric}
\end{proposition}

An operator satisfying the R-RIP satisfies, as a consequence, a number of other properties
when composed with other linear operators defined by the atomic decomposition.

\begin{definition}
Given a set $\Psi = \{\psi_1, \ldots, \psi_{|\Psi|}\} \subset \mathbb{C}^{m \times n}$,
define a linear operator $\mathcal{L}_\Psi: \mathbb{C}^{|\Psi|} \rightarrow \mathbb{C}^{m \times n}$ by
\begin{equation}
\mathcal{L}_\Psi \alpha = \sum_{k=1}^{|\Psi|} \alpha_k \psi_k, \quad \forall \alpha \in \mathbb{C}^{|\Psi|}.
\label{eq:def_L_Psi}
\end{equation}
\end{definition}
It follows from (\ref{eq:def_L_Psi}) that the adjoint operator $\mathcal{L}_\Psi^*: \mathbb{C}^{m \times n} \rightarrow \mathbb{C}^{|\Psi|}$ is given by
\begin{equation}
\left( \mathcal{L}_\Psi^* X \right)_k = \langle X, \psi_k \rangle_{\mathbb{C}^{m \times n}}, \quad \forall k = 1,\ldots, |\Psi|, ~ \forall X \in \mathbb{C}^{m \times n}.
\label{eq:def_adj_L_Psi}
\end{equation}
Note that for $\mathcal{A} : \mathbb{C}^{m \times n} \rightarrow \mathbb{C}^p$
the operator composition $\mathcal{A} \mathcal{L}_\Psi : \mathbb{C}^{|\Psi|} \rightarrow \mathbb{C}^p$ admits a matrix representation.
Its pseudo-inverse is denoted by $\left[ \mathcal{A} \mathcal{L}_\Psi \right]^\dagger$.

\begin{remark}
If $\Psi$ is an orthonormal set, then $\mathcal{L}_\Psi$ is an isometry and $\mathcal{P}_\Psi = \mathcal{L}_\Psi \mathcal{L}_\Psi^*$.
If $\Psi$ is a set of atoms in $\mathbb{O}$, then $\rank(\mathcal{L}_\Psi \alpha) \leq |\Psi|$ for all $\alpha \in \mathbb{C}^{|\Psi|}$.
\end{remark}

\begin{proposition}
Suppose that linear operator $\mathcal{A}: \mathbb{C}^{m \times n} \rightarrow \mathbb{C}^p$ has the rank-restricted isometry constant $\delta_r(\mathcal{A})$.
Let $\Psi$ be a set of atoms in $\mathbb{O}$ such that $|\Psi| \leq r$.
Then
\begin{eqnarray}
\norm{\mathcal{P}_\Psi \mathcal{A}^* b}_F \leq \sqrt{1 + \delta_r(\mathcal{A})} \norm{b}_2, \quad \forall b \in \mathbb{C}^{p}. \label{eq:prop_rip_eq1}
\end{eqnarray}
\label{prop:rip}
\end{proposition}

\begin{proposition}
Suppose that linear operator $\mathcal{A}: \mathbb{C}^{m \times n} \rightarrow \mathbb{C}^p$ has the rank-restricted isometry constant $\delta_r(\mathcal{A})$.
Let $\Psi$ be a set of atoms in $\mathbb{O}$ such that $|\Psi| \leq r$ and let $X \in \mathbb{C}^{m \times n}$ satisfy $\rank(X) \leq r$.
Then
\begin{equation}
\norm{\mathcal{P}_\Psi \mathcal{A}^* \mathcal{A} X}_F \leq (1 + \delta_r(\mathcal{A})) \norm{X}_F.
\end{equation}
\label{prop:rip_proj}
\end{proposition}

\begin{proposition}
Suppose that linear operator $\mathcal{A}: \mathbb{C}^{m \times n} \rightarrow \mathbb{C}^p$ has the rank-restricted isometry constant $\delta_r(\mathcal{A})$.
Let $\Psi$ be a set of atoms in $\mathbb{O}$ such that $|\Psi| \leq r$
and let $\mathcal{P} : \mathbb{C}^{m \times n} \rightarrow \mathbb{C}^{m \times n}$ be a projection operator that commutes with $\mathcal{P}_\Psi$.
Then
\begin{equation}
(1 - \delta_r(\mathcal{A})) \norm{\mathcal{P} \mathcal{P}_\Psi X}_F
\leq \norm{\mathcal{P} \mathcal{P}_\Psi \mathcal{A}^* \mathcal{A} \mathcal{P} \mathcal{P}_\Psi X}_F, \quad \forall X \in \mathbb{C}^{m \times n}.
\end{equation}
\label{prop:rip_proj2}
\end{proposition}

The following \emph{rank-restricted orthogonality property} for the matrix case is analogous to
the sparsity-restricted orthogonality property for the vector case (Lemma~2.1 in \cite{candes2008rip}).

\begin{proposition}
Suppose that linear operator $\mathcal{A}: \mathbb{C}^{m \times n} \rightarrow \mathbb{C}^p$ has the rank-restricted isometry constant $\delta_r(\mathcal{A})$.
Let $X, Y \in \mathbb{C}^{m \times n}$ satisfy $\langle X, Y \rangle_{\mathbb{C}^{m \times n}} = 0$ and $\rank(X + \alpha Y) \leq r$ for all $\alpha \in \mathbb{C}$.
Then
\begin{equation}
\left| \langle \mathcal{A} X , \mathcal{A} Y \rangle_{\mathbb{C}^p} \right| \leq \sqrt{2} \delta_r(\mathcal{A}) \norm{X}_F \norm{Y}_F.
\end{equation}
\label{prop:rop_mat}
\end{proposition}

\begin{remark}
For the vector case, the representation of a vector $x \in \mathbb{C}^n$ in terms of the standard basis $\{e_j\}_{j=1}^n$ of $\mathbb{C}^n$ determines $\norm{x}_0$.
Let $J_1, J_2 \subset \{1,\ldots,n\}$ be arbitrary.
Then the following properties hold:
(i) the projection operators $\mathcal{P}_{\{e_j\}_{j \in J_1}}$ and $\mathcal{P}_{\{e_j\}_{j \in J_2}}$ commute;
and (ii) $\mathcal{P}_{\{e_j\}_{j \in J_1}}^\perp x$ is $s$-sparse (or sparser) if $x$ is $s$-sparse.
These properties follow from the orthogonality of the standard basis.
Proposition 3.2 in \cite{needell2009cis}, corresponding in the vector case to our Proposition~\ref{prop:rop_mat}, requires these two properties.
However, the analogues of properties (i) and (ii) do not hold for the matrix case.
Indeed, for $\Psi_1, \Psi_2 \subset \mathbb{O}$, the projection operators $\mathcal{P}_{\Psi_1}$ and $\mathcal{P}_{\Psi_2}$ do not commute in general
and $\rank(\mathcal{P}_\Psi X)$ can be greater than $r$ even though $\rank(X) \leq r$.
Proposition~\ref{prop:rop_mat} is a stronger version of the corresponding result (Proposition 3.2 in \cite{needell2009cis}) for the vector case
in the sense that it requires a weaker condition (orthogonality between two low-rank matrices),
which can be satisfied without the analogues of properties (i) and (ii).
\end{remark}

\begin{corollary}
Suppose that linear operator $\mathcal{A}: \mathbb{C}^{m \times n} \rightarrow \mathbb{C}^p$ has the rank-restricted isometry constant $\delta_r(\mathcal{A})$.
If sets $\Psi,\Upsilon$ of atoms in $\mathbb{O}$ and matrix $X \in \mathbb{C}^{m \times n}$ satisfy
$\mathcal{P}_\Upsilon^\perp \mathcal{P}_\Psi = \mathcal{P}_\Psi \mathcal{P}_\Upsilon^\perp$, $\mathcal{P}_\Upsilon^\perp X = 0$, and $|\Psi| \leq r$,
then
\begin{equation}
\norm{ \mathcal{P}_\Upsilon^\perp \mathcal{P}_\Psi \mathcal{A}^* \mathcal{A} \mathcal{P}_\Psi X }_F \leq \sqrt{2} \delta_r(\mathcal{A}) \norm{\mathcal{P}_\Psi X}_F.
\end{equation}
\label{cor:rop_perp_proj}
\end{corollary}

\begin{remark}
For the real matrix case, Proposition~\ref{prop:rop_mat} can be improved by dropping the constant $\sqrt{2}$.
This improvement is achieved by replacing the parallelogram identity in the proof to the version for the real scalar field case.
This argument also applies to Corollary~\ref{cor:rop_perp_proj}.
\end{remark}

Finally, we relate the R-RIP to the nuclear norm, extending the analogous result \cite{needell2009cis}
from the $r$-sparse vector case to the rank-$r$ matrix case.
\begin{proposition}
If a linear map $\mathcal{A}: \mathbb{C}^{m \times n} \rightarrow \mathbb{C}^p$ satisfies
\begin{equation}
\norm{\mathcal{A}X}_2^2 \leq (1 + \delta_r(\mathcal{A})) \norm{X}_F^2,
\end{equation}
for all $X \in \mathbb{C}^{m \times n}$ with $\rank(X) \leq r$, then
\begin{equation}
\norm{\mathcal{A}X}_2 \leq \sqrt{1 + \delta_r(\mathcal{A})} \left[ \norm{X}_F + \frac{1}{\sqrt{r}} \norm{X}_* \right],
\end{equation}
for all $X \in \mathbb{C}^{m \times n}$.
\label{prop:energy_bnd}
\end{proposition}

\section{Proof of Theorem~\ref{thm:pg_gen}}
\label{sec:proof_performance_guarantee}

\subsection{Exactly Low Rank Matrix Case}

\begin{theorem}
Assume $\rank(X_0) \leq r$ in (\ref{eq:noisy_meas_mdl}).
Let $\widehat{X}_k$ denote the estimate of $X_0$ in the $k$-th iteration of ADMiRA.
Then for each $k \geq 0$, $\widehat{X}_k$ satisfies the following recursion:
\begin{equation*}
\|X_0 - \widehat{X}_{k+1}\|_F \leq 0.5 \|X_0 - \widehat{X}_k\|_F + 6.5 \norm{\nu}_2.
\end{equation*}
From the above relation, it follows that
\begin{equation*}
\|X_0 - \widehat{X}_k\|_F \leq 2^{-k} \norm{X_0}_F + 13 \norm{\nu}_2, \quad \forall k \geq 0.
\end{equation*}
\label{thm:pg_exact_low_rank}
\end{theorem}

Theorem~\ref{thm:pg_exact_low_rank} is proved by applying a sequence of lemmata.
We generalize the proof of the performance guarantee for CoSaMP \cite{needell2009cis} to the matrix case
by applying the generalized analogy proposed in this paper.
The flow and the techniques used in the proofs are similar to those in \cite{needell2009cis}.
However, in the matrix case, there are additional unknowns in the form of the singular vectors.
Therefore, the generalization of the proofs in \cite{needell2009cis} to the matrix case is not straightforward
and the proofs are sufficiently different from those for the vector case to warrant detailed exposition.
The main steps in the derivation of the performance guarantee are stated in this section and the detailed proofs are in the Appendix.

For the proof, we study the $(k+1)$-th iteration starting with the previous result in the $k$-th iteration.
Let $X_0$ denote the true solution with rank $r$.
Matrix $\widehat{X}$ denotes $\widehat{X}_k$, which is the estimate of $X_0$ in the $k$-th (previous) iteration.
Set $\widehat{\Psi}$ is the set of orthogonal atoms obtained in the previous iteration.
From $(b - \mathcal{A}\widehat{X})$, we compute the proxy matrix $\mathcal{A}^* (b - \mathcal{A}\widehat{X})$.
Set $\Psi'$ is the solution of the following low rank approximation problem:
\begin{equation*}
\Psi' \triangleq \displaystyle \arg\max_\Psi \left\{ \norm{\mathcal{P}_\Psi \mathcal{A}^* (b - \mathcal{A}\widehat{X})}_F :~ \Psi \subset \mathbb{O}, \quad |\Psi| \leq 2r \right\}
\end{equation*}

\begin{lemma}
Let $\rank(X_0) \leq r$ in (\ref{eq:noisy_meas_mdl}). Then
\begin{equation*}
\norm{\mathcal{P}_{\Psi'}^\perp (X_0 - \widehat{X})}_F \leq 0.24 \|X_0 - \widehat{X}\|_F + 2.13 \norm{\nu}_2
\end{equation*}
\label{lemma:identification}
\end{lemma}

Lemma~\ref{lemma:identification} shows that subject to the rank-restricted isometry property,
the set $\Psi'$ of atoms chosen in Step~\ref{step:cor_max_A} of ADMiRA is a good set:
it captures $94\% (= 1 - 0.24^2)$ of the energy of the atoms in $X_0$ that were not captured by $\widehat{X}$,
and the effects of additive measurement noise are bounded by a small constant.
In other words, the algorithm is guaranteed to make good progress in this step.

\begin{lemma}
Let $X_0, \widehat{X} \in \mathbb{C}^{m \times n}$ and
let $\Psi', \widehat{\Psi}$ be sets of atoms in $\mathbb{O}$ such that
$|\Psi'| \leq 2r$, $|\widehat{\Psi}| \leq r$, and $\mathcal{P}_{\widehat{\Psi}}^\perp \widehat{X} = 0$.
Let $\widetilde{\Psi} = \Psi' \cup \widehat{\Psi}$.
Then
\begin{equation*}
\norm{\mathcal{P}_{\widetilde{\Psi}}^\perp X_0}_F \leq \norm{\mathcal{P}_{\Psi'}^\perp (X_0 - \widehat{X})}_F.
\end{equation*}
\label{lemma:merge}
\end{lemma}

Lemma~\ref{lemma:merge} shows that the augmented set of atoms $\widetilde{\Psi}$ produced in Step~\ref{step:merge} of the algorithm is at least as good in explaining the unknown $X_0$
as was the set $\Psi'$ in explaining the part of $X_0$ not captured by the estimate $\widehat{X}$ from the previous iteration.

\begin{lemma}
Let $\rank(X_0) \leq r$ in (\ref{eq:noisy_meas_mdl}) and
let $\widetilde{\Psi}$ be a set of atoms in $\mathbb{O}$ with $|\widetilde{\Psi}| \leq 3r$.
Then
\begin{equation}
\widetilde{X} = \displaystyle \arg\min_{X} \left\{ \norm{b - \mathcal{A}X}_2 :~ X \in \mathrm{span}(\widetilde{\Psi}) \right\}
\label{eq:lemma_estimation}
\end{equation}
satisfies
\begin{equation*}
\|X_0 - \widetilde{X}\|_F \leq 1.04 \norm{\mathcal{P}_{\widetilde{\Psi}}^\perp X_0}_F + 1.02 \norm{\nu}_2.
\end{equation*}
\label{lemma:estimation}
\end{lemma}

Lemma~\ref{lemma:estimation} shows that the least-squares step, Step~\ref{step:solve_LS} of the algorithm, performs almost as well as one could do with operator $\mathcal{A}$ equal to an identity operator:
because $\widetilde{X}$ is restricted to $\mathrm{span}(\widetilde{\Psi})$, it is impossible to recover components of $X_0$ in $\widetilde{\Psi}^\perp$.
Hence, the first constant cannot be smaller than 1.
A value of 1 for the second constant, the noise gain, would correspond to a perfectly conditioned system.

\begin{lemma}
Let $\rank(X_0) \leq r$ in (\ref{eq:noisy_meas_mdl}) and let $\widetilde{X}_r$ denote the best rank-$r$ approximation of $\widetilde{X}$, \textit{i.e.},
\begin{equation*}
\widetilde{X}_r = \arg\min_X \left\{ \|\widetilde{X} - X\|_F :~ \rank(X) \leq r \right\}.
\end{equation*}
Then
\begin{equation*}
\|X_0 - \widetilde{X}_r\|_F \leq 2 \|X_0 - \widetilde{X}\|_F.
\end{equation*}
\label{lemma:prune}
\end{lemma}

As expected, reducing the rank of the estimate $\widetilde{X}$ from $3r$ to $r$,
to produce $\widetilde{X}_r$, increases the approximation error.
However, Lemma~\ref{lemma:prune} shows that this increase is moderate -- by no more than a factor of 2.

The update $\widehat{X}_{k+1} = \widetilde{X}_r$ completes the $(k+1)$-th iteration.
Combining all the results in the lemmata provides the proof of Theorem~\ref{thm:pg_exact_low_rank}.
\begin{IEEEproof} \textbf{(Theorem~\ref{thm:pg_exact_low_rank})}
\begin{eqnarray*}
\|X_0 - \widehat{X}_{k+1}\|_F
&=& \|X_0 - \widetilde{X}_r\|_F \\
&\leq& 2 \|X_0 - \widetilde{X}\|_F \quad \text{(Lemma~\ref{lemma:prune})}\\
&\leq& 2 \cdot \left( 1.04 \norm{\mathcal{P}_{\widetilde{\Psi}}^\perp X_0}_F + 1.02 \norm{\nu}_2 \right) \quad \text{(Lemma~\ref{lemma:estimation})}\\
&\leq& 2.08 \norm{\mathcal{P}_{\Psi'}^\perp (X_0 - \widehat{X}_k)}_F + 2.04 \norm{\nu}_2 \quad \text{(Lemma~\ref{lemma:merge})}\\
&\leq& 2.08 \cdot \left( 0.24 \|X_0 - \widehat{X}_k\|_F + 2.13 \norm{\nu}_2 \right) + 2.04 \norm{\nu}_2 \quad \text{(Lemma~\ref{lemma:identification})}\\
&\leq& 0.5 \|X_0 - \widehat{X}_k\|_F + 6.5 \norm{\nu}_2.
\end{eqnarray*}
The recursion together with the fact that $ \sum_{j=0}^k 2^{-j} \leq \sum_{j=0}^\infty 2^{-j} = 2$ provide the final result.
%\qed
\end{IEEEproof}

\subsection{General Matrix Case}

Theorem~\ref{thm:pg_gen} is proved by combining Theorem~\ref{thm:pg_exact_low_rank} and the following lemma,
which shows how to convert the mismodeling error (deviations of $X_0$ from a low rank matrix) to an equivalent additive measurement noise with a quantified norm.
\begin{lemma}
Let $X_0$ be an arbitrary matrix in $\mathbb{C}^{m \times n}$.
The measurement $b = \mathcal{A}X_0 + \nu$ is also represented as $b = \mathcal{A}X_{0,r} + \widetilde{\nu}$ where
\begin{equation*}
\norm{\widetilde{\nu}}_2 \leq 1.02 \left[ \norm{X_0 - X_{0,r}}_F + \frac{1}{\sqrt{r}} \norm{X_0 - X_{0,r}}_* \right] + \norm{\nu}_2
\end{equation*}
\label{lemma:noise}
\end{lemma}
\begin{IEEEproof}
Let $\widetilde{\nu} = \mathcal{A}(X_0 - X_{0,r}) + \nu$. Then $b = \mathcal{A}X_{0,r} + \widetilde{\nu}$.
\begin{eqnarray*}
\norm{\widetilde{\nu}}_2
&\leq& \norm{\mathcal{A}(X_0 - X_{0,r})}_2 + \norm{\nu}_2 \\
&\leq& \sqrt{1 + \delta_r(\mathcal{A})} \left[ \norm{X_0 - X_{0,r}}_F + \frac{1}{\sqrt{r}} \norm{X_0 - X_{0,r}}_* \right] + \norm{\nu}_2,
\end{eqnarray*}
where the last inequality holds by Proposition~\ref{prop:energy_bnd}.
The inequality $\delta_r(\mathcal{A}) \leq \delta_{4r}(\mathcal{A}) \leq 0.04$ implies $\sqrt{1 + \delta_r(\mathcal{A})} \leq 1.02$.
%\qed
\end{IEEEproof}

\begin{IEEEproof} \textbf{(Theorem~\ref{thm:pg_gen})}
Let $X$ be an arbitrary matrix in $\mathbb{C}^{m \times n}$.
The measurement is given by $b = \mathcal{A}X_{0,r} + \widetilde{\nu}$, where $\widetilde{\nu}$ is defined in Lemma~\ref{lemma:noise}.
By Theorem~\ref{thm:pg_exact_low_rank},
\begin{equation*}
\|X_{0,r} - \widehat{X}_{k+1}\|_F \leq 0.5 \|X_{0,r} - \widehat{X}_k\|_F + 6.5 \norm{\widetilde{\nu}}_2.
\end{equation*}
Applying the triangle inequality and the above inequality,
\begin{eqnarray*}
\|X_0 - \widehat{X}_{k+1}\|_F
&\leq& \|X_{0,r} - \widehat{X}_{k+1}\|_F + \norm{X_0 - X_{0,r}}_F \\
&\leq& 0.5 \|X_{0,r} - \widehat{X}_k\|_F + 6.5 \norm{\widetilde{\nu}}_2 + \norm{X_0 - X_{0,r}}_F \\
\end{eqnarray*}
Using the upper bound on $\norm{\widetilde{\nu}}_2$ yields
\begin{eqnarray*}
\|X_0 - \widehat{X}_{k+1}\|_F
&\leq& 0.5 \|X_0 - \widehat{X}_k\|_F + 7.63 \norm{X_0 - X_{0,r}}_F + \frac{6.63}{\sqrt{r}} \norm{X_0 - X_{0,r}}_* + 6.5 \norm{\nu}_2 \\
&<& 0.5 \|X_0 - \widehat{X}_k\|_F + 8 \epsilon,
\end{eqnarray*}
where $\epsilon$ is the unrecoverable energy.
%\qed
\end{IEEEproof}

\section{Required Number of Iterations}
\label{sec:iter_cnt}

Theorem~\ref{thm:iter_cnt_gen} provides a uniform bound on the number of iterations required to achieve the guaranteed approximation accuracy.
In addition to this uniform iteration bound,
Theorem~\ref{thm:bnd_iter_gen} in this section shows that even faster convergence may be expected for matrices $X$ with clustered singular values.

In the analysis of the iteration number, the distribution of the singular values of the matrices involved is the only thing that matters.
Indeed, the singular vectors do not play any role in the analysis.
As a consequence, the proofs for the vector case (CoSaMP) and the matrix case (ADMiRA) are very similar, and the corresponding bounds on the number of iterations coincide.
However, for the completeness, we provide the proofs for the matrix case.

\begin{definition}
Given $X \in \mathbb{C}^{m \times n}$, $\atoms{X}$ is defined in (\ref{eq:def_atom}).
We define the atomic bands of $X$ by
\begin{equation*}
B_j \triangleq \{ \psi \in \atoms{X} :~ 2^{-(j+1)} \norm{X}_F^2 < \norm{\mathcal{P}_\psi X}_F^2 \leq 2^{-j} \norm{X}_F^2 \}, \quad \textrm{for} ~ j \in \mathbb{Z}_+,
\end{equation*}
where $\mathbb{Z}_+$ denotes the set of nonnegative integers.
Note that atomic bands are disjoint subsets of $\atoms{X}$, which is an orthonormal set of atoms in $\mathbb{O}$, and therefore atomic bands are mutually orthogonal.
From the atomic bands, the profile of $X$ is defined as the number of nonempty atomic bands, \textit{i.e.},
\begin{equation}
\mathrm{profile}(X) \triangleq |\{ j :~ B_j \neq \emptyset \}|.
\label{eq:def_profile}
\end{equation}
\end{definition}
From the definition, $\mathrm{profile}(X) \leq \rank(X)$.

The atomic bands and $\mathrm{profile}(X)$ admit a simple interpretation in terms of the spectrum of $X$.
Let $\atoms{X} = \{ \psi_k \}_{k=1}^{\rank(X)}$ be ordered as $\norm{\mathcal{P}_{\psi_k} X}_F \geq \norm{\mathcal{P}_{\psi_{k+1}} X}_F$.
Then $\norm{\mathcal{P}_{\Psi_k} X}_F = \sigma_k$, where $\sigma_k$ is the $k$-th singular value of $X$, in decreasing order.
Let
\begin{equation*}
\widetilde{\sigma}_k^2 = \frac{\sigma_k^2}{\norm{X}_F^2}.
\end{equation*}
Then $B_j = \{ \psi_k :~ -(j+1) \leq \log_2 \widetilde{\sigma}_k^2 \leq -j \}$.
In other words, $B_j$ contains the $\atoms{X}$ corresponding to normalized singular values falling in a one octave interval (``bin'').
The quantity $\mathrm{profile}(X)$ then is the number of such occupied octave bins, and measures the spread of singular values of $X$ on a log scale.

\begin{remark}
For the vector case, the term analogous to the atomic band is the component band \cite{needell2008cisTR} defined by
\begin{equation*}
B_j \triangleq \{ k \in \{1,\ldots,n\} :~ 2^{-(j+1)} \norm{x}_2^2 < |x_k|^2 \leq 2^{-j} \norm{x}_2^2 \}, \quad \textrm{for} ~ j \in \mathbb{Z}_+,
\end{equation*}
for $x \in \mathbb{C}^n$.
\end{remark}

First, the number of iterations for the exactly low-rank case is bounded by the following theorem.
\begin{theorem}
Let $X_0 \in \mathbb{C}^{m \times n}$ be a rank-$r$ matrix and let $t = \mathrm{profile}(X_0)$,
where $\mathrm{profile}(X_0)$ is defined in (\ref{eq:def_profile}).
Then after at most
\begin{equation*}
t \log_{4/3}(1 + 4.3 \sqrt{r/t}) + 6
\end{equation*}
iterations, the estimate $\widehat{X}$ produced by ADMiRA satisfies
\begin{equation*}
\|X_0 - \widehat{X}\|_F \leq 15 \norm{\nu}_2.
\end{equation*}
\label{thm:bnd_iter_exact}
\end{theorem}

We introduce additional notations for the proof of Theorem~\ref{thm:bnd_iter_exact}.
Let $\widehat{X}_k$ denote the estimate at the $k$-th iteration of ADMiRA.
For a nonnegative integer $j$, we define an auxiliary matrix
\begin{equation*}
Y_j \triangleq \sum_{k \geq j} \mathcal{P}_{B_k} X_0.
\end{equation*}
Then $Y_j$ satisfies
\begin{equation}
\norm{Y_j}_F^2 \leq \sum_{k \geq j} 2^{-k} \norm{X_0}_F^2 \cdot |B_k|.
\label{eq:bnd_norm_Yj}
\end{equation}

The proof of Theorem~\ref{thm:bnd_iter_exact} is done by a sequence of lemmata.
The first lemma presents two possibilities in each iteration of ADMiRA:
if the iteration is successful, the approximation error is small;
otherwise, the approximation error is dominated by the un-identified portion of the matrix
and the approximation error in the next iteration decreases by a constant ratio.

\begin{lemma}
Let $\rank(X_0) \leq r$ in (\ref{eq:noisy_meas_mdl}).
Matrix $\widehat{X}_k$ denotes the estimate of $X_0$ in the $k$-th iteration of ADMiRA.
Let $\widehat{\Psi}_k$ denote $\atoms{\widehat{X}_k}$.
In each iteration of ADMiRA, at least one of the followings holds:
either
\begin{equation}
\|X_0 - \widehat{X}_k\|_F \leq 70 \norm{\nu}_2,
\label{eq:first_alt}
\end{equation}
or
\begin{eqnarray}
\|X_0 - \widehat{X}_k\|_F &\leq& 2.15 \norm{ \mathcal{P}_{\widehat{\Psi}_k}^\perp X }_F
\label{eq:second_alt_a}
\\
\|X_0 - \widehat{X}_{k+1}\|_F &\leq& \left(\frac{3}{4}\right) \|X_0 - \widehat{X}_k\|_F.
\label{eq:second_alt_b}
\end{eqnarray}
\label{lemma:two_possibilities}
\end{lemma}

\begin{lemma}
Fix $K = \lfloor t \log_{4/3}(1 + 4.3\sqrt{r/t}) \rfloor$.
Assume that (\ref{eq:second_alt_a}) and (\ref{eq:second_alt_b}) are in force for each iteration.
Then $\atoms{\widehat{X}_K} = \atoms{X_0}$.
\label{lemma:atom_recovery}
\end{lemma}

Next, the result is extended to the approximately low-rank case by using, once again, Lemma~\ref{lemma:noise}.
\begin{theorem}
Let $X_0 \in \mathbb{C}^{m \times n}$ be an arbitrary matrix and let $t = \mathrm{profile}(X_{0,r})$.
Then, after at most
\begin{equation*}
t \log_{4/3}(1 + 4.3 \sqrt{r/t}) + 6
\end{equation*}
iterations, the estimate $\widehat{X}$ produced by ADMiRA satisfies
\begin{equation*}
\|X_0 - \widehat{X}\|_F \leq 17 \epsilon,
\end{equation*}
where $\epsilon$ is the unrecoverable energy.
\label{thm:bnd_iter_gen}
\end{theorem}

\begin{IEEEproof} \textbf{(Theorem~\ref{thm:iter_cnt_gen})}
As a function of $t$, $(t \log_{4/3}(1 + 4.3 \sqrt{r/t}) + 6)$ is maximized when $t = r$.
Since $\log_{4/3}5.6 < 6$, the number of iterations is at most $6(r+1)$.
Therefore, the approximation error of ADMiRA is achieved within $6(r+1)$ iterations for any matrix $X_0$.
%\qed
\end{IEEEproof}

Theorem~\ref{thm:bnd_iter_gen} is also of independent interest, because the bound it provides reveals that
even faster convergence can be achieved for matrices $X$ with small $\mathrm{profile}(X_{0,r}) \ll r$.
Recall the relationship between $\mathrm{profile}(X_{0,r})$ and the distribution of the $r$ largest singular values of $X_0$.
It follows that the number of iterations in ADMiRA required for convergence is roughly proportional to the number of clusters of singular values of $X_{0,r}$ on a log scale.

\section{Implementation and Scalability}
\label{sec:large}

We analyze the computational complexity of ADMiRA and will show that ADMiRA scales well to large problem instances.
Each iteration of ADMiRA consists of procedures requiring the following basic operations:
application of $\mathcal{A}$ and $\mathcal{A}^*$, singular value decompositions, and solving a least square problem.
We analyze the computational cost of the procedures in terms of the complexity of the basic operations, which will depend on the properties of $\mathcal{A}$.
First note that ADMiRA keeps the matrix variables (except the proxy matrix) in factorized form through their atomic decomposition,
which is advantageous for both the computational efficiency and memory requirements.
Furthermore, the proxy matrix is often sparse in applications such as the matrix completion problem.

\textit{Computing the proxy matrix:}
this involves the application of $\mathcal{A}$ and $\mathcal{A}^*$.
The procedure first computes the residual $y = b - \mathcal{A} \widehat{X}$ and then computes the proxy matrix $\mathcal{A}^* y$.
Let $\widehat{X} = \sum_{k=1}^r \sigma_k u_k v_k^H$ denote the atomic decomposition of $\widehat{X}$.
Here $u_k v_k^H$'s are not necessarily orthogonal.
$(\mathcal{A} \widehat{X})_k$ can be computed by
$\langle \widehat{X}, Z_k \rangle_{\mathbb{C}^{m \times n}} = \sum_{k=1}^r \sigma_k v_k^H Z_k^H u_k$, $k=1,\ldots,p$,
for an appropriate set of $p$ matrices $Z_k \in \mathbb{C}^{m \times n}$.
Then $\mathcal{A}^* y$ can be computed by $\sum_{k=1}^p y_k Z_k$.
The complexity of these operations will depend on the sparsity of $\mathcal{A}$.
\begin{description}
\item[Case 1]: $\mathcal{A}$ is an arbitrary linear (dense) operator
and the costs of computing $\mathcal{A} \widehat{X}$ and $\mathcal{A}^* y$ are $O(prmn)$ and $O(pmn)$, respectively.
\item[Case 2]: $\mathcal{A}$ is a sparse linear operator -- so the $Z_k$ have $O(m+n)$ non-zero elements, and
and the costs of computing $\mathcal{A} \widehat{X}$ and $\mathcal{A}^* y$ are $O(pr(m+n))$ and $O(p(m+n))$, respectively.
\item[Case 3]: $\mathcal{A}$ is an extremely sparse linear operator (such as in the matrix completion problem),
so the $Z_k$ have $O(1)$ nonzeros,
and the costs of computing $\mathcal{A} \widehat{X}$ and $\mathcal{A}^* y$ are $O(pr)$ and $O(p)$, respectively.
\end{description}

\textit{Finding the $2r$ principal atoms of the proxy matrix:}
this involves the truncated singular value decomposition with $2r$ dominant singular triplets, which can be computed by the Lanczos method
at a cost of $O(mnrL)$, where $L$ denotes the number of the Lanczos iterations per each singular value, which depends on the singular value distribution.
An alternative approach is to use recent advances in low rank approximation of large matrices based on randomized algorithms
(c.f. \cite{harpeled2006lrm}, \cite{woolfe2007fra}, and the references therein.)
that compute the low-rank approximation of a given matrix in time linear in the size of the matrix.
These randomized algorithms are useful when the size of the matrix is large but the rank $r$ remains a small constant.
For example, the complexity of Har-Peled's algorithm \cite{harpeled2006lrm} is $O(mnr^2\log r)$.
When $\mathcal{A}$ is sparse with $O(1)$ nonzero elements per each $Z_k$,
the matrix-vector product $(\mathcal{A}^*y)w$ for $w \in \mathbb{C}^p$ can be computed as
$\sum_{k=1}^p y_k Z_k w$ and hence the complexity reduces to $O(prL)$ for the Lanczos method and $O(pr^2 \log r)$ for the randomized method, respectively.

\textit{Solving least square problems:}
ADMiRA requires the solution of an over-determined system with $p$ equations and $3r$ unknowns.
The complexity is $O(pr^2)$. Similarly to CoSaMP, the Richardson iteration or the conjugate gradient method can be used to improve the complexity of this part.
The convergence of the Richardson iteration is guaranteed owing to the R-RIP assumption of ADMiRA and the complexity is $O(pr)$.

\textit{Finding the $r$ principal atoms of the solution to the least square problem:}
this also involves the truncated singular value decomposition of the least square solution $\widetilde{X}$.
In fact, this procedure can be done more efficiently by exploiting the fact that
$\widetilde{X}$ is available in a factorized form $\widetilde{X} = U \Sigma V^H$ where
$U \in \mathbb{C}^{m \times 3r}, V \in \mathbb{C}^{n \times 3r}$, and $\Sigma$ is a $3r \times 3r$ diagonal matrix.
Here $U, V$ do not consist of orthogonal columns in general.
Let $U = Q_U R_U$ and $V = Q_V R_V$ denote the QR factorizations of $U$ and $V$, respectively.
Then $Q_U^H Q_U = I_m$ and $Q_V^H Q_V = I_n$.
Now let $W D Z^H$ denote the singular value decomposition of the $3r \times 3r$ matrix $R_U \Sigma R_V^H$.
Then we have the desired singular value decomposition $\widetilde{X} = (Q_U W) D (Q_V Z)^H$.
The complexity is $O((m+n+r)r^2)$, which is negligible compared to a direct SVD of $\widetilde{X}$.

Applications of $\mathcal{A}$ and $\mathcal{A}^*$ are the most demanding procedures of ADMiRA for a dense linear operator $\mathcal{A}$.
These operations are also required in all other algorithms for $\text{P1}$, $\text{P1'}$, or $\text{P2}$.
To overcome this computational complexity, the linear operator $\mathcal{A}$ should have some structure that admits efficient computation.
Examples include random Toeplitz matrices and randomly subsampled Fourier measurements.
For matrix completion, $\mathcal{A}$ is a sparse with $O(1)$ cost per measurement
and hence these operations are dominated by the remaining operations.
In this case, the computation of the truncated singular value decomposition is the most demanding procedure of ADMiRA.
Equipped with the randomized low rank approximation,
ADMiRA has complexity of $O(p r^2\log r)$ per iteration, or $O(p r^3\log r)$ to achieve the guarantee in Theorem~\ref{thm:iter_cnt_gen}.
ADMiRA therefore has complexity linear in the size $p$ of the data, and it scales well to large problems.

\section{Numerical Experiment}
\label{sec:num_exp}

\begin{figure}[htb]
\begin{center}
\begin{minipage}[htb]{0.45\linewidth}
\centerline{\includegraphics[height=50mm]{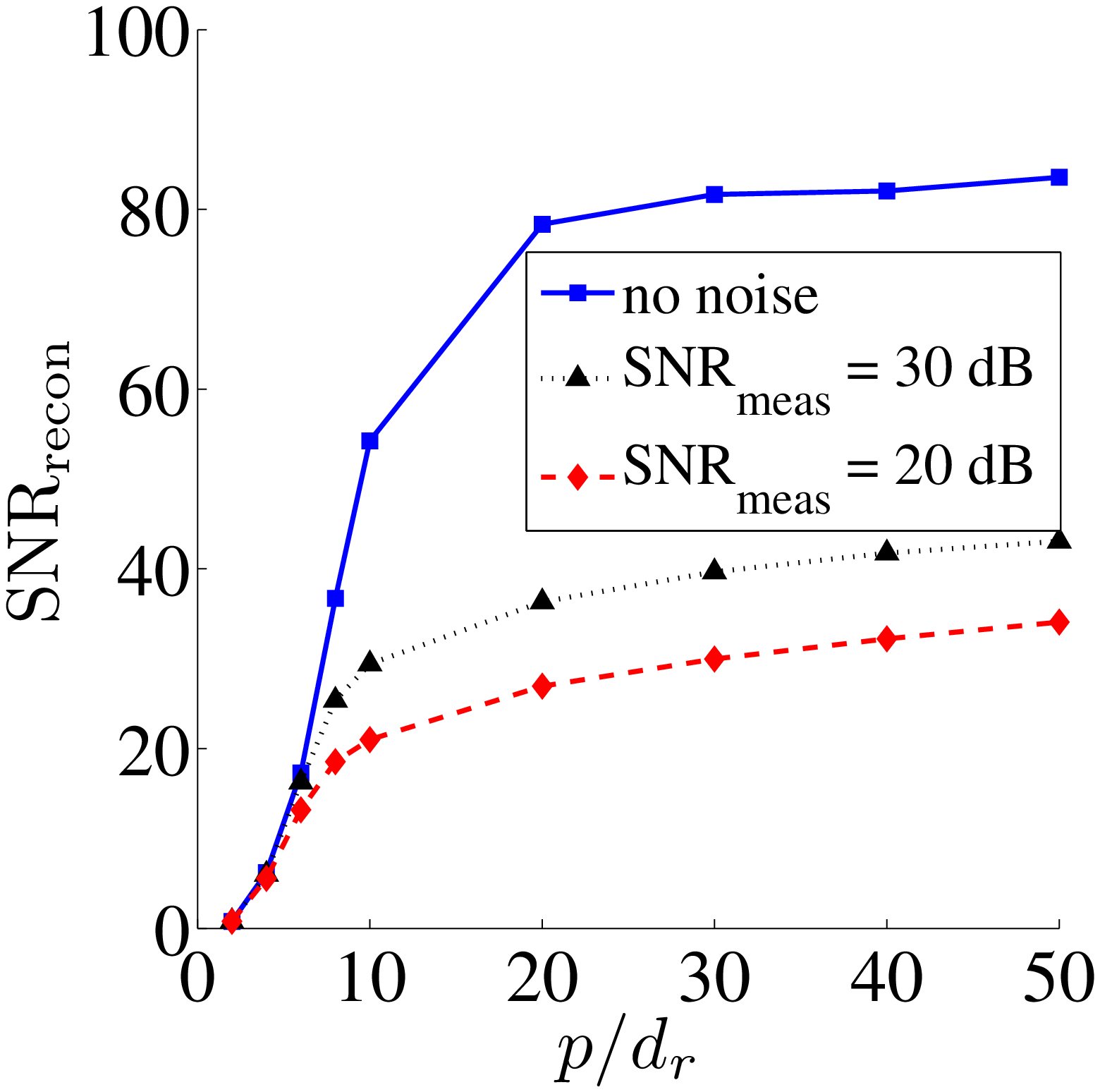}}
\end{minipage}
\begin{minipage}[htb]{0.45\linewidth}
\centerline{\includegraphics[height=50mm]{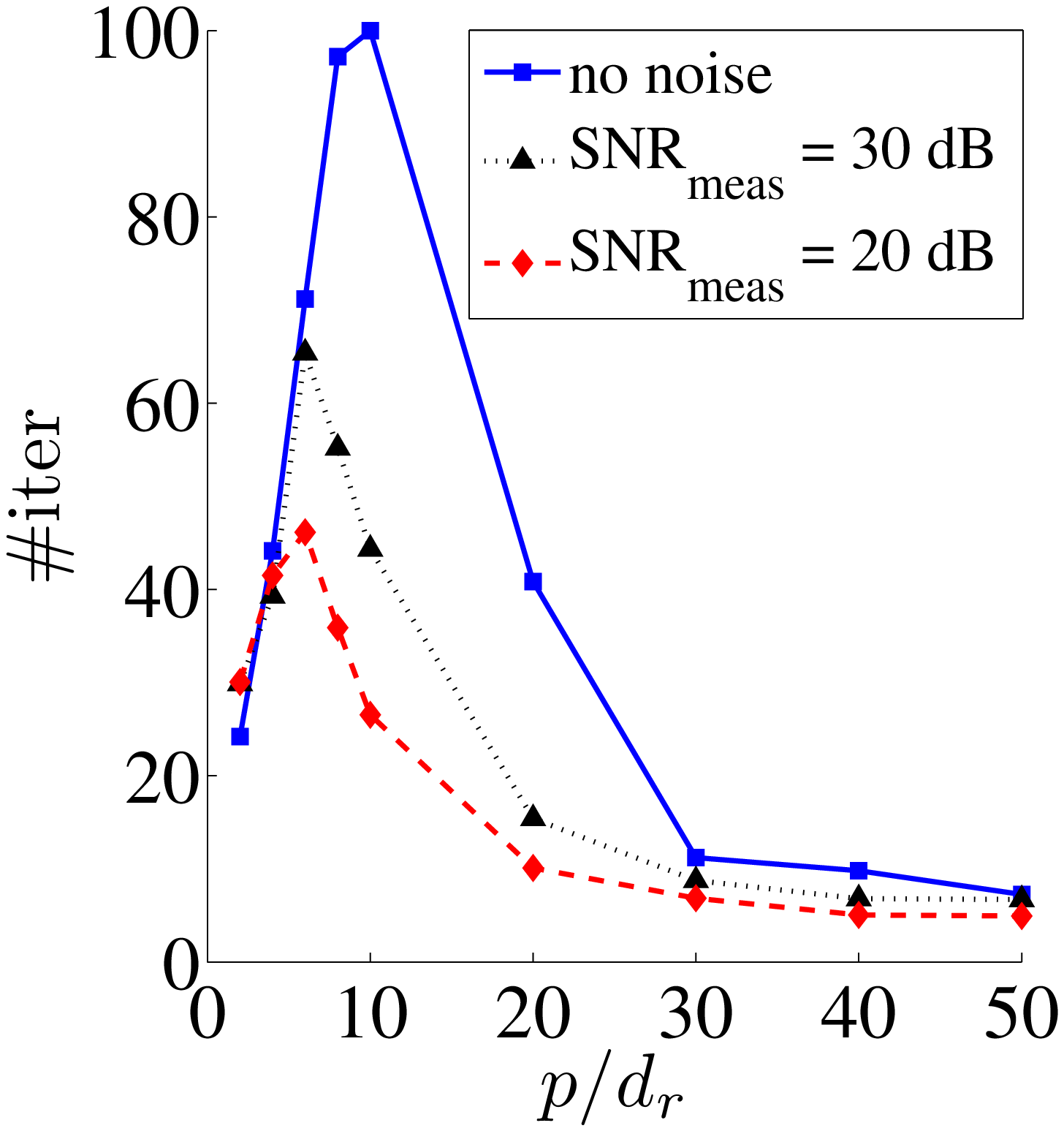}}
\end{minipage}
\end{center}
\caption{Completion of random matrices by ADMiRA: $n = m = 500, r = 2$.}
\label{fig:var_p}
\end{figure}

We tested the performance of ADMiRA with an operator $\mathcal{A}$ generated by a Gaussian ensemble,
which satisfies RIP with high probability.
ADMiRA performed well in this case as predicted by our theory.
Here we study reconstructions by ADMiRA with a generic matrix completion example.
Note that the performance guarantee in terms of R-RIP does not applies to this case,
because the linear operator in the matrix completion problem does not satisfy the RIP.
None the less, we want to check the empirical performance of ADMiRA in this practically important application.
Our Matlab implementation uses PROPACK \cite{larsen:psp} (an implementation of the Lanczos algorithm) to compute partial SVDs in Steps~\ref{step:cor_max_A} and \ref{step:cor_max_C} of ADMiRA.
The test matrix $X_0 \in \mathbb{R}^{n \times n}$ is generated as the product $X_0 = Y_L Y_R^H$
where $Y_L,Y_R \in \mathbb{R}^{n \times r}$ have entries following an i.i.d. Gaussian distribution.
The measurement $b$ is $p$ randomly chosen entries of $X$, which may be contaminated with an additive white Gaussian noise.
The reconstruction error and measurement noise level are measured in terms of
$\mathrm{SNR}_\mathrm{recon} \triangleq 20 \log_{10}(\norm{X_0}_F / \|X_0 - \widehat{X}\|_F)$ and
$\mathrm{SNR}_\mathrm{meas} \triangleq 20 \log_{10}(\norm{b}_2 / \norm{\nu}_2)$, respectively.
Computational efficiency is measured by the number of iterations.
Here we stopped the algorithm when $\|b - \mathcal{A} \widehat{X}\|_2 / \norm{b}_2 < 10^{-4}$.
As a result, the algorithm provided $\mathrm{SNR}_\mathrm{recon}$ around $70 \mathrm{dB}$ for the ideal (noiseless and exactly low-rank) case when it was successful.
However, it is still possible to get higher $\mathrm{SNR}_\mathrm{recon}$ with a few more iterations.
The results in Fig.~\ref{fig:var_p}, Table~\ref{table:var_n}, and Table~\ref{table:compare_admira_svt} have been averaged over 20 trials.

Fig.~\ref{fig:var_p} shows that both $\mathrm{SNR}_\mathrm{recon}$ and the number of iterations improve as $p/d_r$ increases.
Here $d_r$ is the number of degrees of freedom in a real rank-$r$ matrix defined by $d_r = r(n+m-r)$ and denotes the essential number of unknowns.
Fig.~\ref{fig:var_p} suggests that we need $p / d_r \geq 20$ for $n = 500$.

Candes and Recht \cite{candes2008emc} showed that $p = O(n^{1.2} r \log_{10} n)$ known entries suffice to complete an unknown $n \times n$ rank-$r$ matrix.
Table~\ref{table:var_n} shows that ADMiRA provides nearly perfect recovery of random matrices from $p$ known entries where $p = 10 \lceil n^{1.2} r \log_{10} n \rceil$.
Although $\mathrm{SNR}_\mathrm{recon}$ in the noiseless measurement case is high enough to say that the completion is nearly perfect,
the number of iterations increases as $n$ increases.
We are studying whether this increase in iterations with $n$ might be an artifact of our numerical implementation of ADMiRA.
In the noisy measurement case the number of iterations is low and does not increase with problem size $n$.
Because in most if not all practical applications the data will be noisy, or the matrix to be recovered only approximately low rank,
this low and constant number of iterations is of practical significance.

Table~\ref{table:compare_admira_svt} shows that in most of the examples tested, ADMiRA provides slightly better performance with less computation than SVT \cite{cai2008svt}.
Roughly, the computational complexity of a single iteration of ADMiRA can be compared to two times that of SVT.

\begin{figure}[htb]
\begin{center}
\begin{minipage}[htb]{0.45\linewidth}
\centerline{\includegraphics[height=50mm]{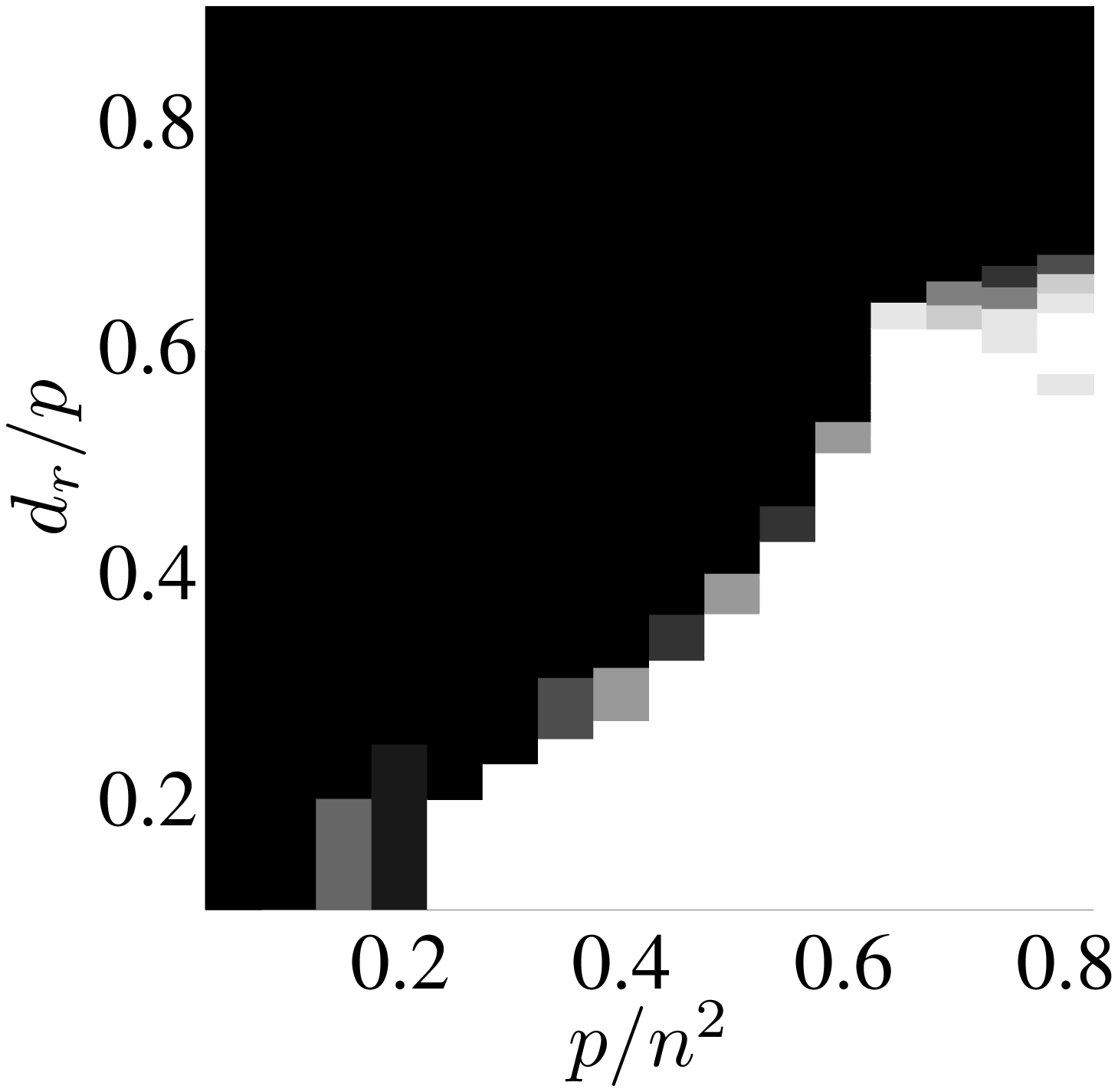}}
\centering{\footnotesize ADMiRA}
\end{minipage}
\begin{minipage}[htb]{0.45\linewidth}
\centerline{\includegraphics[height=50mm]{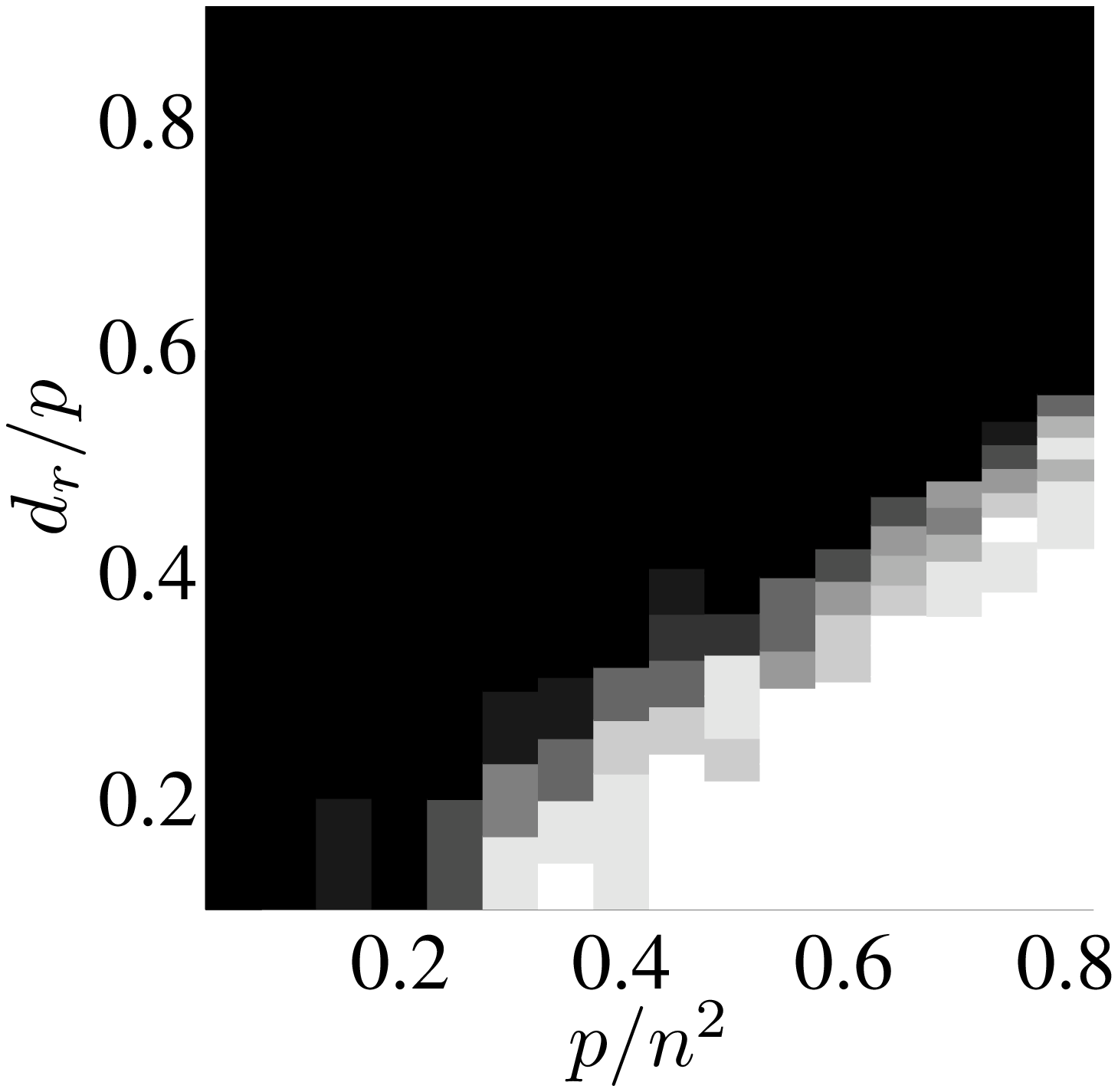}}
\centering{\footnotesize SVT}
\end{minipage}
\end{center}
\caption{Phase transition of matrix completion: $n = m = 100$.}
\label{fig:ptd}
\end{figure}

Fig.~\ref{fig:ptd} compares the phase transitions of ADMiRA and SVT.
We count the number of successful matrix completions ($\mathrm{SNR}_\mathrm{recon} \geq 70 \mathrm{dB}$) out of 10 trials for each triplet $(n,p,r)$.
Brighter color implies more success.
ADMiRA performed better than SVT for this example.

We emphasize that all comparisons with SVT were performed for the noiseless exactly low rank matrix case,
because the current implementation \cite{svt:software} and theory \cite{cai2008svt} of SVT do not support the ellipsoidal constraint case.
We are not aware of an efficient, scalable algorithm other than ADMiRA that supports the ellipsoidal constraint.

\begin{table}
\begin{center}
\begin{minipage}[b]{0.9\linewidth}
\centerline{\begin{tabular}{c|*{5}{c|}c}
  \hline
  \multirow{2}{*}{$n$} & \multirow{2}{*}{$p/n^2$} & \multirow{2}{*}{$p/d_r$} & \multicolumn{2}{|c|}{no noise} & \multicolumn{2}{|c}{$\mathrm{SNR}_\mathrm{meas} = 20 \mathrm{dB}$} \\\cline{4-7}
  & & & $\mathrm{SNR}_\mathrm{recon}$ (dB) & $\#\text{iter}$ & $\mathrm{SNR}_\mathrm{recon}$ (dB) & $\#\text{iter}$ \\\hline\hline
  500 & 0.37 & 47 & 83 & 8 & 34 & 5 \\\hline
  1000 & 0.24 & 60 & 83 & 9 & 34 & 5 \\\hline
  1500 & 0.18 & 69 & 82 & 11 & 35 & 5 \\\hline
  2000 & 0.15 & 76 & 81 & 12 & 35 & 5 \\\hline
  2500 & 0.13 & 81 & 81 & 18 & 36 & 5 \\\hline
  3000 & 0.12 & 86 & 81 & 24 & 36 & 5 \\\hline
  3500 & 0.10 & 90 & 81 & 26 & 36 & 5 \\\hline
  4000 & 0.09 & 95 & 80 & 32 & 36 & 5 \\\hline
  4500 & 0.09 & 98 & 81 & 37 & 36 & 5 \\\hline
\end{tabular}}
\end{minipage}
\end{center}
\caption{Completion of random matrices by ADMiRA: $n = m$, $r = 2$, $p = 10 \lceil n^{1.2} r \log_{10} n \rceil$.}
\label{table:var_n}
\end{table}

\begin{table}
\begin{center}
\begin{minipage}[b]{0.9\linewidth}
\centerline{\begin{tabular}{c|*{5}{c|}c}
  \hline
  \multirow{2}{*}{$r$} & \multirow{2}{*}{$p/n^2$} & \multirow{2}{*}{$p/d_r$} & \multicolumn{2}{|c|}{$\mathrm{SNR}_\mathrm{recon}$ (dB)} & \multicolumn{2}{|c}{$\#\text{iter}$} \\\cline{4-7}
  & & & ADMiRA & SVT & ADMiRA & SVT \\\hline\hline
  \multirow{6}{*}{2} & 0.05 & 12.51 & 77 & 74 & 259 & 143 \\\cline{2-7}
  & 0.10 & 25.03 & ~~~~79~~~~ & ~~~~77~~~~ & ~~~~56~~~~ & ~~~~77~~~~ \\\cline{2-7}
  & 0.15 & 37.54 & 81 & 78 & 20 & 61 \\\cline{2-7}
  & 0.20 & 50.05 & 82 & 79 & 11 & 54 \\\cline{2-7}
  & 0.25 & 62.56 & 84 & 79 & 8 & 49 \\\cline{2-7}
  & 0.30 & 75.08 & 84 & 79 & 7 & 46 \\\hline
  \multirow{6}{*}{5} & 0.05 & 5.01 & 19 & 37 & 99 & 500 \\\cline{2-7}
  & 0.10 & 10.03 & 77 & 76 & 89 & 100 \\\cline{2-7}
  & 0.15 & 15.04 & 78 & 77 & 32 & 75 \\\cline{2-7}
  & 0.20 & 20.05 & 81 & 78 & 15 & 64 \\\cline{2-7}
  & 0.25 & 25.06 & 82 & 79 & 11 & 57 \\\cline{2-7}
  & 0.30 & 30.08 & 83 & 79 & 8 & 53 \\\hline
  \multirow{6}{*}{10} & 0.05 & 2.51 & 7 & -9 & 28 & 451 \\\cline{2-7}
  & 0.10 & 5.03 & 30 & 74 & 194 & 205 \\\cline{2-7}
  & 0.15 & 7.54 & 77 & 76 & 50 & 99 \\\cline{2-7}
  & 0.20 & 10.05 & 79 & 77 & 19 & 80 \\\cline{2-7}
  & 0.25 & 12.56 & 80 & 78 & 13 & 69 \\\cline{2-7}
  & 0.30 & 15.08 & 80 & 78 & 10 & 62 \\\hline
\end{tabular}}
\end{minipage}
\end{center}
\caption{Comparison of ADMiRA and SVT: no noise, $n = m = 1000$.}
\label{table:compare_admira_svt}
\end{table}

\section{Conclusion}
\label{sec:conclusion}

We proposed a new algorithm, ADMiRA,
which extends both the efficiency and the performance guarantee of the CoSaMP algorithm for $\ell_0$-norm minimization to matrix rank minimization.
The proposed generalized correlation maximization can be also applied to MP, OMP, and SP and their variants
to similarly extend the known algorithms and theory from the $s$-term vector approximation problem to the rank-$r$ matrix approximation.
ADMiRA can handle large scale rank minimization problems efficiently by using recent linear time algorithms for low rank approximation of a known matrix..
Our numerical experiments demonstrate that ADMiRA is an effective algorithm even when the R-RIP is not satisfied, as in the matrix completion problem.
While the performance guarantee in this paper relies on the R-RIP, it seems that a performance guarantee for ADMiRA without using the R-RIP might be possible.

\appendix

\subsection{Proof of Proposition~\ref{prop:monotone_ric}}
The rank-restricted isometry constant $\delta_r(\mathcal{A})$ can be represented as
\begin{equation}
\delta_r(\mathcal{A}) = \max\{ [\sigma_{r,\max}(\mathcal{A})]^2 - 1 , 1 - [\sigma_{r,\min}(\mathcal{A})]^2 \},
\label{eq:ric_comp}
\end{equation}
where $\sigma_{r,\max}(\mathcal{A})$ and $\sigma_{r,\min}(\mathcal{A})$ are defined by
\begin{equation*}
\sigma_{r,\max}(\mathcal{A}) \triangleq \max_X \{ \norm{\mathcal{A}X}_2 :~ \norm{X}_F = 1, \quad \rank(X) \leq r \},
\end{equation*}
and
\begin{equation*}
\sigma_{r,\min}(\mathcal{A}) \triangleq \min_X \{ \norm{\mathcal{A}X}_2 :~ \norm{X}_F = 1, \quad \rank(X) \leq r \},
\end{equation*}
respectively.
As $r$ increases, the feasible sets of both problems increase
and hence $\sigma_{r,\max}(\mathcal{A})$ and $\sigma_{r,\min}(\mathcal{A})$ are nondecreasing and nonincreasing, respectively.
Therefore, (\ref{eq:ric_comp}) implies that $\delta_r(\mathcal{A})$ is nondecreasing in $r$.

\subsection{Proof of Proposition~\ref{prop:rip}}
Let $d = \dim(\mathrm{span}(\Psi))$ and let $\Phi = \{\phi_j\}_{j=1}^d$ be an orthonormal basis of $\mathrm{span}(\Psi)$.
\footnote{
Note that $\Phi$ is not necessarily a set of atoms in $\mathbb{O}$.
}
Then $\mathcal{L}_\Phi$ is an isometry that satisfies $\mathcal{P}_\Psi = \mathcal{L}_\Psi \mathcal{L}_\Psi^\dagger = \mathcal{L}_\Phi \mathcal{L}_\Phi^*$.

Since $\rank(\mathcal{L}_\Phi \alpha) \leq |\Psi| \leq r$ and $\norm{\mathcal{L}_\Phi \alpha}_F = \norm{\alpha}_2$ for all $\alpha \in \mathbb{C}^d$, by the R-RIP
\begin{equation}
\norm{\mathcal{A} \mathcal{L}_\Phi \alpha}_2
\leq \sqrt{1 + \delta_r(\mathcal{A})} \norm{\mathcal{L}_\Phi \alpha}_F
= \sqrt{1 + \delta_r(\mathcal{A})} \norm{\alpha}_2, \quad \forall \alpha \in \mathbb{C}^d.
\label{eq:proof_prop_rip}
\end{equation}
This implies that the operator norm of $\mathcal{A} \mathcal{L}_\Phi$ is bounded from above by $\sqrt{1 + \delta_r(\mathcal{A})}$.
Since the adjoint operator $\left[\mathcal{A} \mathcal{L}_\Phi\right]^*$ has the same operator norm,
\begin{eqnarray}
\norm{\mathcal{L}_\Phi \left[\mathcal{A} \mathcal{L}_\Phi\right]^* b}_F = \norm{\left[\mathcal{A} \mathcal{L}_\Phi\right]^* b}_2 \leq \sqrt{1 + \delta_r(\mathcal{A})} \norm{b}_2, \quad \forall b \in \mathbb{C}^p. \label{eq:proof_prop_rip_eq1}
\end{eqnarray}
Then (\ref{eq:prop_rip_eq1}) follows from (\ref{eq:proof_prop_rip_eq1}) with $\mathcal{L}_\Phi \left[\mathcal{A} \mathcal{L}_\Phi\right]^* = \mathcal{L}_\Phi \mathcal{L}_\Phi^* \mathcal{A}^* = \mathcal{P}_\Psi \mathcal{A}^*$.

\subsection{Proof of Proposition~\ref{prop:rip_proj}}
Let $Y = \mathcal{P}_\Psi \mathcal{A}^* \mathcal{A} X$. Then $\rank(Y) \leq |\Psi| \leq r$.
By R-RIP,
\begin{eqnarray*}
\left| \langle \mathcal{A} X, \mathcal{A} Y \rangle_{\mathbb{C}^p} \right|^2
\leq \norm{\mathcal{A} X}_2^2 \norm{\mathcal{A} Y}_2^2
\leq (1 + \delta_r(\mathcal{A}))^2 \norm{X}_F^2 \norm{Y}_F^2.
\end{eqnarray*}
Therefore
\begin{eqnarray*}
\langle \mathcal{A} X, \mathcal{A} Y \rangle_{\mathbb{C}^p}
= \langle \mathcal{A} X, \mathcal{A} \mathcal{P}_\Psi \mathcal{A}^* \mathcal{A} X \rangle_{\mathbb{C}^p}
= \langle \mathcal{P}_\Psi \mathcal{A}^* \mathcal{A} X, \mathcal{P}_\Psi \mathcal{A}^* \mathcal{A} X \rangle_{\mathbb{C}^{m \times n}}
= \norm{\mathcal{P}_\Psi \mathcal{A}^* \mathcal{A} X}_F^2
\leq (1 + \delta_r(\mathcal{A})) \norm{X}_F \norm{\mathcal{P}_\Psi \mathcal{A}^* \mathcal{A} X}_F.
\end{eqnarray*}

\subsection{Proof of Proposition~\ref{prop:rip_proj2}}
Since $P$ and $\mathcal{P}_\Psi$ are commuting projection operators,
$P \mathcal{P}_\Psi$ is a projection operator onto $S = \mathcal{R}(P) \cap \mathcal{R}(\mathcal{P}_\Psi)$.
Let $d = \dim(S)$ and let $\Phi = \{\phi_k\}_{k=1}^d \subset \mathbb{C}^{m \times n}$ be an orthonormal basis of $S$,
then $\mathcal{L}_\Phi$ is an isometry that satisfies $P \mathcal{P}_\Psi = \mathcal{L}_\Phi \mathcal{L}_\Phi^*$.
Since $S \subset \mathcal{R}(\mathcal{P}_\Psi)$, $\rank(\mathcal{L}_\Phi \alpha) \leq |\Psi| \leq r$ for all $\alpha \in \mathbb{C}^d$.
Therefore, by R-RIP,
\begin{equation*}
\sqrt{1 - \delta_r(\mathcal{A})} \norm{\alpha}_2 = \sqrt{1 - \delta_r(\mathcal{A})} \norm{\mathcal{L}_\Phi \alpha}_F \leq \norm{\mathcal{A} \mathcal{L}_\Phi \alpha}_2, \quad \forall \alpha \in \mathbb{C}^d,
\end{equation*}
where the first equality holds since $\mathcal{L}_\Phi$ is an isometry.
By the relationship between $\mathcal{A} \mathcal{L}_\Phi$ and $\mathcal{L}_\Phi^* \mathcal{A}^* \mathcal{A} \mathcal{L}_\Phi$, it follow that
\begin{equation*}
(1 - \delta_r(\mathcal{A})) \norm{\alpha}_2 \leq \norm{\mathcal{L}_\Phi^* \mathcal{A}^* \mathcal{A} \mathcal{L}_\Phi \alpha}_2, \quad \forall \alpha \in \mathbb{C}^d.
\end{equation*}
For each $X \in \mathbb{C}^{m \times n}$, there exists $\alpha \in \mathbb{C}^d$ such that $P \mathcal{P}_\Psi X = \mathcal{L}_\Phi \alpha$.
\begin{eqnarray*}
\norm{\mathcal{P} \mathcal{P}_\Psi \mathcal{A}^* \mathcal{A} \mathcal{P} \mathcal{P}_\Psi X}_F
&=& \norm{\mathcal{L}_\Phi \mathcal{L}_\Phi^* \mathcal{A}^* \mathcal{A} \mathcal{L}_\Phi \alpha}_F
= \norm{\mathcal{L}_\Phi^* \mathcal{A}^* \mathcal{A} \mathcal{L}_\Phi \alpha}_2 \\
&\geq& (1 - \delta_r(\mathcal{A})) \norm{\alpha}_2
= (1 - \delta_r(\mathcal{A})) \norm{\mathcal{L}_\Phi \alpha}_F
= (1 - \delta_r(\mathcal{A})) \norm{\mathcal{P} \mathcal{P}_\Psi X}_F.
\end{eqnarray*}

\subsection{Proof of Proposition~\ref{prop:rop_mat}}
Assume that $\norm{X}_F = \norm{Y}_F = 1$.
Let $\alpha \in \mathbb{C}$ be a constant of unit modulus, \textit{i.e.} $|\alpha| = 1$.
By the subadditivity of the rank, $\rank(X + \alpha Y) \leq r$.
By the orthogonality of $X$ and $Y$, $\norm{X + \alpha Y}_F^2 = \norm{X}_F^2 + |\alpha| \norm{Y}_F^2 = 2$.
Therefore
\begin{equation*}
2(1 - \delta_r(\mathcal{A})) \leq \norm{\mathcal{A} X + \alpha \mathcal{A} Y}_2^2 \leq 2(1 + \delta_r(\mathcal{A})).
\end{equation*}
In particular, the inequality holds for $\alpha = \pm 1, \pm i$ where $i = \sqrt{-1}$.
By the parallelogram identity,
\begin{eqnarray*}
\left| \langle \mathcal{A} X , \mathcal{A} Y \rangle_{\mathbb{C}^p} \right|^2
&=& \frac{1}{16} \left|
\norm{\mathcal{A} X + \mathcal{A} Y}_2^2 - \norm{\mathcal{A} X - \mathcal{A} Y}_2^2
+ i \norm{\mathcal{A} X + i \mathcal{A} Y}_2^2 - i \norm{\mathcal{A} X - i \mathcal{A} Y}_2^2
\right|^2 \\
&=& \frac{1}{16} \left| \norm{\mathcal{A} X + \mathcal{A} Y}_2^2 - \norm{\mathcal{A} X - \mathcal{A} Y}_2^2 \right|^2
+ \frac{1}{16} \left| \norm{\mathcal{A} X + i \mathcal{A} Y}_2^2 - \norm{\mathcal{A} X - i \mathcal{A} Y}_2^2 \right|^2 \\
&\leq& 2 \left[ \delta_r(\mathcal{A}) \right]^2.
\end{eqnarray*}

\subsection{Proof of Corollary~\ref{cor:rop_perp_proj}}
For an arbitrary matrix $Y \in \mathbb{C}^{m \times n}$,
\begin{equation*}
\langle \mathcal{P}_\Psi X , \mathcal{P}_\Psi \mathcal{P}_\Upsilon^\perp Y \rangle_{\mathbb{C}^{m \times n}}
= \langle X , \mathcal{P}_\Psi \mathcal{P}_\Upsilon^\perp Y \rangle_{\mathbb{C}^{m \times n}}
= \langle X , \mathcal{P}_\Upsilon^\perp \mathcal{P}_\Psi Y \rangle_{\mathbb{C}^{m \times n}}
= \langle \mathcal{P}_\Upsilon^\perp X , \mathcal{P}_\Upsilon^\perp \mathcal{P}_\Psi Y \rangle_{\mathbb{C}^{m \times n}} = 0
\end{equation*}
and
\begin{equation*}
\rank(\mathcal{P}_\Psi X + \alpha \mathcal{P}_\Psi \mathcal{P}_\Upsilon^\perp Y) \leq \rank(\mathcal{P}_\Psi (X + \alpha \mathcal{P}_\Upsilon^\perp Y)) \leq |\Psi| \leq r
\end{equation*}
for all $\alpha \in \mathbb{C}$.
Therefore Proposition~\ref{prop:rop_mat} implies
\begin{eqnarray*}
\langle \mathcal{A} \mathcal{P}_\Psi X , \mathcal{A} \mathcal{P}_\Psi \mathcal{P}_\Upsilon^\perp Y \rangle_{\mathbb{C}^p} \leq \sqrt{2} \delta_r(\mathcal{A}) \norm{\mathcal{P}_\Psi X}_F \norm{\mathcal{P}_\Psi \mathcal{P}_\Upsilon^\perp Y}_F.
\end{eqnarray*}
Since $Y$ was arbitrary, we can take $Y = \mathcal{A}^* \mathcal{A} \mathcal{P}_\Psi X$.
Then
\begin{eqnarray*}
\langle \mathcal{A} \mathcal{P}_\Psi X , \mathcal{A} \mathcal{P}_\Psi \mathcal{P}_\Upsilon^\perp Y \rangle_{\mathbb{C}^p}
&=& \langle \mathcal{P}_\Upsilon^\perp \mathcal{P}_\Psi \mathcal{A}^* \mathcal{A} \mathcal{P}_\Psi X , \mathcal{P}_\Psi \mathcal{P}_\Upsilon^\perp \mathcal{A}^* \mathcal{A} \mathcal{P}_\Psi X \rangle_{\mathbb{C}^{m \times n}} \\
&=& \norm{\mathcal{P}_\Upsilon^\perp \mathcal{P}_\Psi \mathcal{A}^* \mathcal{A} \mathcal{P}_\Psi X}_F^2 \\
&\leq& \sqrt{2} \delta_r(\mathcal{A}) \norm{\mathcal{P}_\Psi X}_F \norm{\mathcal{P}_\Upsilon^\perp \mathcal{P}_\Psi \mathcal{A}^* \mathcal{A} \mathcal{P}_\Psi X}_F.
\end{eqnarray*}

\subsection{Proof of Proposition~\ref{prop:energy_bnd}}
We modify the proof the analogous result for the vector case in \cite{needell2009cis} for our proposition.

For $\Psi \subset \mathbb{O}$, the unit-ball in the subspace spanned by $\Psi$ is defined by
\begin{equation*}
B_F^\Psi \triangleq \{ X \in \mathbb{C}^{m \times n} :~ X \in \mathrm{span}(\Psi), \quad \norm{X}_F \leq 1 \}.
\end{equation*}
Define the convex body
\begin{equation*}
S \triangleq \mathrm{conv} \left\{ \bigcup_{\Psi \subset \mathbb{O}, |\Psi| \leq r} B_F^\Psi \right\},
\end{equation*}
where $\mathrm{conv}\{G\}$ denotes the convex hull of set $G$.
By the assumption, the operator norm satisfies
\begin{equation*}
\norm{\mathcal{A}}_{S \rightarrow 2} \triangleq \max_{X \in S} \norm{\mathcal{A}X}_2 \leq \sqrt{1 + \delta_r(\mathcal{A})}.
\end{equation*}
Define the second convex body
\begin{equation*}
K \triangleq \left\{ X \in \mathbb{C}^{m \times n} :~ \norm{X}_F + \frac{1}{\sqrt{r}} \norm{X}_* \leq 1 \right\},
\end{equation*}
and consider the operator norm
\begin{equation*}
\norm{\mathcal{A}}_{K \rightarrow 2} \triangleq \max_{X \in K} \norm{\mathcal{A}X}_2.
\end{equation*}
The claim of the proposition is equivalent to
\begin{equation*}
\norm{\mathcal{A}}_{K \rightarrow 2} \leq \norm{\mathcal{A}}_{S \rightarrow 2}.
\end{equation*}
It suffices to show that $K \subset S$.
Let $X$ be an element in $K$.
Consider the singular value decomposition of $X$,
\begin{equation*}
X = \sum_{k=1}^{\rank(X)} \sigma_k u_k v_k^H,
\end{equation*}
with $\sigma_{k+1} \leq \sigma_k$.
Let $\sigma_k = 0$ if $k > \rank(X)$ and $J = \lceil \rank(X) / r \rceil - 1$,
where $\lceil c \rceil$ is the smallest integer equal to or greater than $c$.
Then we have the following decomposition
\begin{equation*}
X = \sum_{j=0}^J \sum_{k=rj+1}^{r(j+1)} \sigma_k u_k v_k^H = \sum_{j=0}^J c_j Y_j,
\end{equation*}
where
\begin{equation*}
c_j \triangleq \norm{\sum_{k=rj+1}^{r(j+1)} \sigma_k u_k v_k^H}_F \quad \mathrm{and} \quad
Y_j \triangleq \frac{1}{c_j} \sum_{k=rj+1}^{r(j+1)} \sigma_k u_k v_k^H.
\end{equation*}

For each $j \in \{1,\ldots,J\}$,
\begin{eqnarray*}
c_j &=& \sqrt{\sum_{k=rj+1}^{r(j+1)} \sigma_k^2} \leq \sqrt{r} \cdot \sigma_{rj+1} \leq \sqrt{r} \cdot \frac{1}{r} \sum_{k=r(j-1)+1}^{rj} \sigma_k.
\end{eqnarray*}
Therefore
\begin{eqnarray*}
\sum_{j=1}^J c_j \leq \frac{1}{\sqrt{r}} \sum_{j=1}^J \sum_{k=r(j-1)+1}^{rj} \sigma_k = \frac{1}{\sqrt{r}} \sum_{k=1}^{\rank(X)} \sigma_k = \frac{1}{\sqrt{r}} \norm{X}_*.
\end{eqnarray*}
From the definition of $c_0$, it follows that $c_0 \leq \norm{X}_F$.
Since $X \in K$, we note
\begin{eqnarray*}
\sum_{j=0}^J c_j \leq \norm{X}_F + \frac{1}{\sqrt{r}} \norm{X}_* \leq 1.
\end{eqnarray*}

Also note that $Y_j \in S$ for all $j = 0,\ldots,J$ since $\rank(Y_j) \leq r$ and $\norm{Y_j}_F = 1$ by construction.
Therefore $X$ is the convex combination of the elements in $S$.
Since $S$ is a convex hull, $X \in S$.

\subsection{Proof of Lemma~\ref{lemma:identification}}
Let $\Phi = \atoms{X_0 - \widehat{X}}$.
Since $|\Phi| \leq \rank(X_0) + \rank(\widehat{X}) \leq 2r$,
it follows by the selection rule of $\Psi'$ that
\begin{equation}
\norm{\mathcal{P}_\Phi \mathcal{A}^* (b - \mathcal{A}\widehat{X})}_F \leq \norm{\mathcal{P}_{\Psi'} \mathcal{A}^* (b - \mathcal{A}\widehat{X})}_F.
\label{eq:ineq1_lemma_identification}
\end{equation}
Let $\Upsilon \subset \mathbb{O}$ be a set of atoms that spans $\mathrm{span}(\Phi) \cap \mathrm{span}(\Psi')$.
Then $\mathcal{P}_\Upsilon$ and $\mathcal{P}_\Phi$ commute,
\begin{equation}
\mathcal{P}_\Upsilon \mathcal{P}_\Phi = \mathcal{P}_\Phi \mathcal{P}_\Upsilon \quad \mathrm{and} \quad \mathcal{P}_\Upsilon^\perp \mathcal{P}_\Phi = \mathcal{P}_\Phi \mathcal{P}_\Upsilon^\perp,
\end{equation}
as do $\mathcal{P}_\Upsilon$ and $\mathcal{P}_{\Psi'}$,
\begin{equation}
\mathcal{P}_\Upsilon \mathcal{P}_{\Psi'} = \mathcal{P}_{\Psi'} \mathcal{P}_\Upsilon \quad \mathrm{and} \quad \mathcal{P}_\Upsilon^\perp \mathcal{P}_{\Psi'} = \mathcal{P}_{\Psi'} \mathcal{P}_\Upsilon^\perp.
\end{equation}
From the commutativity,
we note that $\mathcal{P}_\Upsilon^\perp \mathcal{P}_\Phi$ and $\mathcal{P}_\Upsilon^\perp \mathcal{P}_{\Psi'}$ are projection operators.
Furthermore, each of the projection operators $\mathcal{P}_\Phi$ and $\mathcal{P}_{\Psi'}$ can be decomposed as the sum of two mutually orthogonal projection operators:
\begin{eqnarray}
\mathcal{P}_\Phi = \mathcal{P}_\Upsilon + \mathcal{P}_\Upsilon^\perp \mathcal{P}_\Phi
\label{eq:decomp_proj_Phi}
\\
\mathcal{P}_\Psi' = \mathcal{P}_\Upsilon + \mathcal{P}_\Upsilon^\perp \mathcal{P}_{\Psi'}.
\label{eq:decomp_proj_Psi'}
\end{eqnarray}

Applying (\ref{eq:decomp_proj_Phi}) and (\ref{eq:decomp_proj_Psi'}) to (\ref{eq:ineq1_lemma_identification}), invoking the Pythagorean theorem,
and removing the common term containing $\mathcal{P}_\Upsilon$ gives
\begin{equation}
\norm{\mathcal{P}_\Upsilon^\perp \mathcal{P}_\Phi \mathcal{A}^* (b - \mathcal{A}\widehat{X})}_F \leq \norm{\mathcal{P}_\Upsilon^\perp \mathcal{P}_{\Psi'} \mathcal{A}^* (b - \mathcal{A}\widehat{X})}_F.
\label{eq:ineq2_lemma_identification}
\end{equation}

First, we derive an upper bound on the right hand side of inequality (\ref{eq:ineq2_lemma_identification}).
\begin{eqnarray}
\norm{\mathcal{P}_\Upsilon^\perp \mathcal{P}_{\Psi'} \mathcal{A}^* (b - \mathcal{A}\widehat{X})}_F
&=& \norm{\mathcal{P}_\Upsilon^\perp \mathcal{P}_{\Psi'} \mathcal{A}^* \left( \mathcal{A}(X - \widehat{X}) + \nu \right)}_F \nonumber \\
&\leq& \norm{\mathcal{P}_\Upsilon^\perp \mathcal{P}_{\Psi'} \mathcal{A}^* \mathcal{A}(X - \widehat{X})}_F + \norm{\mathcal{P}_\Upsilon^\perp \mathcal{P}_{\Psi'} \mathcal{A}^* \nu}_F.
\label{eq:ineq2ub_lemma_identification}
\end{eqnarray}
Using Proposition~\ref{prop:rip} and $|\Psi'| \leq 2r$, the second term of (\ref{eq:ineq2ub_lemma_identification}) is bounded by
\begin{eqnarray*}
\norm{\mathcal{P}_\Upsilon^\perp \mathcal{P}_{\Psi'} \mathcal{A}^* \nu}_F
\leq \norm{\mathcal{P}_{\Psi'} \mathcal{A}^* \nu}_F
\leq \sqrt{1 + \delta_{2r}(\mathcal{A})} \norm{\nu}_2.
\end{eqnarray*}
The first term of (\ref{eq:ineq2ub_lemma_identification}) is further bounded by
\begin{eqnarray*}
\norm{\mathcal{P}_\Upsilon^\perp \mathcal{P}_{\Psi'} \mathcal{A}^* \mathcal{A}(X_0 - \widehat{X})}_F
&=& \norm{\mathcal{P}_\Upsilon^\perp \mathcal{P}_{\Psi'} \mathcal{A}^* \mathcal{A} (\mathcal{P}_\Upsilon + \mathcal{P}_\Upsilon^\perp \mathcal{P}_\Phi) (X_0 - \widehat{X})}_F \\
&\leq& \norm{\mathcal{P}_\Upsilon^\perp \mathcal{P}_{\Psi'} \mathcal{A}^* \mathcal{A} \mathcal{P}_\Upsilon (X_0 - \widehat{X})}_F
+ \norm{\mathcal{P}_\Upsilon^\perp \mathcal{P}_{\Psi'} \mathcal{A}^* \mathcal{A} \mathcal{P}_\Phi \mathcal{P}_\Upsilon^\perp (X_0 - \widehat{X})}_F \\
&\leq& \norm{\mathcal{P}_\Upsilon^\perp \mathcal{P}_{\Psi'} \mathcal{A}^* \mathcal{A} \mathcal{P}_{\Psi'} \mathcal{P}_\Upsilon (X_0 - \widehat{X})}_F
+ \norm{\mathcal{P}_{\Psi'} \mathcal{A}^* \mathcal{A} \mathcal{P}_\Phi \mathcal{P}_\Upsilon^\perp (X_0 - \widehat{X})}_F \\
&\leq& \sqrt{2} \delta_{2r}(\mathcal{A}) \norm{\mathcal{P}_\Upsilon (X_0 - \widehat{X})}_F
+ (1 + \delta_{2r}(\mathcal{A})) \norm{\mathcal{P}_\Upsilon^\perp (X_0 - \widehat{X})}_F,
\end{eqnarray*}
where the third inequality follows from Corollary~\ref{cor:rop_perp_proj} with $\mathcal{P}_\Upsilon^\perp \mathcal{P}_{\Psi'} = \mathcal{P}_{\Psi'} \mathcal{P}_\Upsilon^\perp$ and $\mathcal{P}_\Upsilon^\perp \mathcal{P}_\Upsilon (X_0 - \widehat{X}) = 0$
and Proposition~\ref{prop:rip_proj} with $|\Psi'| \leq 2r$ and $\rank(\mathcal{P}_\Phi \mathcal{P}_\Upsilon^\perp (X_0 - \widehat{X})) \leq 2r$.

Combining the previous results, we have the following upper bound on the right hand side of inequality (\ref{eq:ineq2_lemma_identification}).
\begin{eqnarray}
\norm{\mathcal{P}_\Upsilon^\perp \mathcal{P}_{\Psi'} \mathcal{A}^* (b - \mathcal{A}\widehat{X})}_F
&\leq& \sqrt{2} \delta_{2r}(\mathcal{A}) \norm{\mathcal{P}_\Upsilon (X_0 - \widehat{X})}_F
+ (1 + \delta_{2r}(\mathcal{A})) \norm{\mathcal{P}_\Upsilon^\perp (X_0 - \widehat{X})}_F \nonumber \\
&& + \sqrt{1 + \delta_{2r}(\mathcal{A})} \norm{\nu}_2.
\label{eq:ineq2ubmore_lemma_identification}
\end{eqnarray}

Next, we derive a lower bound on the left hand side of inequality (\ref{eq:ineq2_lemma_identification}).
\begin{eqnarray}
\norm{\mathcal{P}_\Upsilon^\perp \mathcal{P}_\Phi \mathcal{A}^* (b - \mathcal{A}\widehat{X})}_F
&=& \norm{\mathcal{P}_\Upsilon^\perp \mathcal{P}_\Phi \mathcal{A}^* \left( \mathcal{A}(X_0 - \widehat{X}) + \nu \right)}_F \nonumber \\
&\geq& \norm{\mathcal{P}_\Upsilon^\perp \mathcal{P}_\Phi \mathcal{A}^* \mathcal{A}(X_0 - \widehat{X})}_F - \norm{\mathcal{P}_\Upsilon^\perp \mathcal{P}_\Phi \mathcal{A}^* \nu}_F.
\label{eq:ineq2lb_lemma_identification}
\end{eqnarray}
Using Proposition~\ref{prop:rip}, the second term of (\ref{eq:ineq2lb_lemma_identification}) is further bounded by
\begin{eqnarray*}
- \norm{\mathcal{P}_\Upsilon^\perp \mathcal{P}_\Phi \mathcal{A}^* \nu}_F
\geq - \norm{\mathcal{P}_\Phi \mathcal{A}^* \nu}_F
\geq - \sqrt{1 + \delta_{2r}(\mathcal{A})} \norm{\nu}_2.
\end{eqnarray*}
The first term of (\ref{eq:ineq2lb_lemma_identification}) is further bounded by
\begin{eqnarray}
\norm{\mathcal{P}_\Upsilon^\perp \mathcal{P}_\Phi \mathcal{A}^* \mathcal{A}(X_0 - \widehat{X})}_F
&=& \norm{\mathcal{P}_\Upsilon^\perp \mathcal{P}_\Phi \mathcal{A}^* \mathcal{A} (\mathcal{P}_\Upsilon + \mathcal{P}_\Upsilon^\perp \mathcal{P}_\Phi) (X_0 - \widehat{X})}_F \nonumber \\
&\geq& \norm{\mathcal{P}_\Upsilon^\perp \mathcal{P}_\Phi \mathcal{A}^* \mathcal{A} \mathcal{P}_\Phi \mathcal{P}_\Upsilon^\perp (X_0 - \widehat{X})}_F
- \norm{\mathcal{P}_\Phi \mathcal{P}_\Upsilon^\perp \mathcal{A}^* \mathcal{A} \mathcal{P}_\Upsilon (X_0 - \widehat{X})}_F \nonumber \\
&=& \norm{\mathcal{P}_\Upsilon^\perp \mathcal{P}_\Phi \mathcal{A}^* \mathcal{A} \mathcal{P}_\Phi \mathcal{P}_\Upsilon^\perp (X_0 - \widehat{X})}_F
- \norm{\mathcal{P}_\Upsilon^\perp \mathcal{P}_\Phi \mathcal{A}^* \mathcal{A} \mathcal{P}_\Phi \mathcal{P}_\Upsilon (X_0 - \widehat{X})}_F \nonumber \\
&\geq& (1 - \delta_{2r}(\mathcal{A})) \norm{\mathcal{P}_\Upsilon^\perp (X_0 - \widehat{X})}_F
- \sqrt{2} \delta_{2r}(\mathcal{A}) \norm{\mathcal{P}_\Upsilon (X_0 - \widehat{X})}_F \label{eq:ineq2lb_intermediate1_lemma_identification} \\
&\geq& (1 - \delta_{2r}(\mathcal{A})) \norm{\mathcal{P}_{\Psi'}^\perp (X_0 - \widehat{X})}_F
- \sqrt{2} \delta_{2r}(\mathcal{A}) \norm{\mathcal{P}_\Upsilon (X_0 - \widehat{X})}_F,
\label{eq:ineq2lb_intermediate2_lemma_identification}
\end{eqnarray}
where the inequality (\ref{eq:ineq2lb_intermediate1_lemma_identification}) follows from
Proposition~\ref{prop:rip_proj2} with $\mathcal{P}_\Upsilon^\perp \mathcal{P}_\Phi = \mathcal{P}_\Phi \mathcal{P}_\Upsilon^\perp$ and $|\Phi| \leq 2r$ and Corollary~\ref{cor:rop_perp_proj},
and the last inequality (\ref{eq:ineq2lb_intermediate2_lemma_identification}) follows from the fact that $(\Psi')^\perp \subset (\Upsilon)^\perp$.

Combining the previous results, we have the following lower bound on the left hand side of inequality (\ref{eq:ineq2_lemma_identification}).
\begin{eqnarray}
\norm{\mathcal{P}_\Upsilon^\perp \mathcal{P}_\Phi \mathcal{A}^* (b - \mathcal{A}\widehat{X})}_F
&\geq& (1 - \delta_{2r}(\mathcal{A})) \norm{\mathcal{P}_{\Psi'}^\perp (X_0 - \widehat{X})}_F
- \sqrt{2} \delta_{2r}(\mathcal{A}) \norm{\mathcal{P}_\Upsilon (X_0 - \widehat{X})}_F \nonumber \\
&& - \sqrt{1 + \delta_{2r}(\mathcal{A})} \norm{\nu}_2.
\label{eq:ineq2lbmore_lemma_identification}
\end{eqnarray}

Combining (\ref{eq:ineq2_lemma_identification}), (\ref{eq:ineq2ubmore_lemma_identification}) and (\ref{eq:ineq2lbmore_lemma_identification}) yields
\begin{eqnarray*}
\norm{\mathcal{P}_{\Psi'}^\perp (X_0 - \widehat{X})}_F
&\leq& \frac{1 + \delta_{2r}(\mathcal{A})}{1 - \delta_{2r}(\mathcal{A})} \norm{\mathcal{P}_\Upsilon^\perp (X_0 - \widehat{X})}_F
+ \frac{2 \sqrt{2} \delta_{2r}(\mathcal{A})}{1 - \delta_{2r}(\mathcal{A})} \norm{\mathcal{P}_\Upsilon (X_0 - \widehat{X})}_F
+ \frac{2\sqrt{1 + \delta_{2r}(\mathcal{A})}}{1 - \delta_{2r}(\mathcal{A})} \norm{\nu}_2 \nonumber \\
&\leq& \frac{1}{1 - \delta_{2r}(\mathcal{A})} \frac{4\sqrt{2}\delta_{2r}(\mathcal{A})(1 + \delta_{2r}(\mathcal{A}))}{\sqrt{(1 + \delta_{2r}(\mathcal{A}))^2 + 8 \delta_{2r}(\mathcal{A})^2}} \|X_0 - \widehat{X}\|_F
+ \frac{2\sqrt{1 + \delta_{2r}(\mathcal{A})}}{1 - \delta_{2r}(\mathcal{A})} \norm{\nu}_2,
\end{eqnarray*}
where the second inequality is obtained by maximizing over $\mathcal{P}_\Upsilon$ with the constraint
\begin{equation*}
\norm{\mathcal{P}_\Upsilon^\perp (X_0 - \widehat{X})}_F^2 + \norm{\mathcal{P}_\Upsilon (X_0 - \widehat{X})}_F^2 = \|X_0 - \widehat{X}\|_F^2.
\end{equation*}

Substituting $\delta_{2r}(\mathcal{A}) \leq \delta_{4r}(\mathcal{A}) \leq 0.04$ gives the constants in the final inequality.

\subsection{Proof of Lemma~\ref{lemma:merge}}
Since $\widehat{\Psi} \subset \widetilde{\Psi}$,
$\mathcal{P}_{\widehat{\Psi}}^\perp \widehat{X} = 0$ implies $\mathcal{P}_{\widetilde{\Psi}}^\perp \widehat{X} = 0$ and hence
\begin{equation*}
\norm{\mathcal{P}_{\widetilde{\Psi}}^\perp X}_F = \norm{\mathcal{P}_{\widetilde{\Psi}}^\perp (X_0 - \widehat{X})}_F \leq \norm{\mathcal{P}_{\Psi'}^\perp (X_0 - \widehat{X})}_F,
\end{equation*}
where the inequality holds since $\Psi' \subset \widetilde{\Psi}$ implies $(\widetilde{\Psi})^\perp \subset (\Psi')^\perp$.

\subsection{Proof of Lemma~\ref{lemma:estimation}}
Assume that $\widetilde{\Psi}$ is a linearly independent set of atoms in $\mathbb{O}$.
Otherwise, we can take $\widetilde{\Psi}$ as a maximal linearly independent subset of $\widetilde{\Psi}$.

The minimizer in (\ref{eq:lemma_estimation}) is given by
\begin{equation*}
\widetilde{X} = \mathcal{L}_{\widetilde{\Psi}} \left[ \mathcal{A} \mathcal{L}_{\widetilde{\Psi}} \right]^\dagger b
= \mathcal{L}_{\widetilde{\Psi}} \left[ \mathcal{A} \mathcal{L}_{\widetilde{\Psi}} \right]^\dagger (\mathcal{A} X_0 + \nu).
\end{equation*}

By the triangle inequality,
\begin{eqnarray}
\|X_0 - \widetilde{X}\|_F
&\leq& \|\underbrace{X_0 - \mathcal{L}_{\widetilde{\Psi}} \left[ \mathcal{A} \mathcal{L}_{\widetilde{\Psi}} \right]^\dagger \mathcal{A} X_0}_{\rank \leq 4r}\|_F
+ \|\underbrace{\mathcal{L}_{\widetilde{\Psi}} \left[ \mathcal{A} \mathcal{L}_{\widetilde{\Psi}} \right]^\dagger \nu}_{\rank \leq 3r}\|_F \\
&\leq& \frac{1}{\sqrt{1 - \delta_{4r}(\mathcal{A})}} \norm{\mathcal{A} \left( X_0 - \mathcal{L}_{\widetilde{\Psi}} \left[ \mathcal{A} \mathcal{L}_{\widetilde{\Psi}} \right]^\dagger \mathcal{A} X_0 \right)}_2
+ \frac{1}{\sqrt{1 - \delta_{3r}(\mathcal{A})}} \norm{\mathcal{A} \mathcal{L}_{\widetilde{\Psi}} \left[ \mathcal{A} \mathcal{L}_{\widetilde{\Psi}} \right]^\dagger \nu}_2,
\label{eq:lemma_estimation_proof_eq1}
\end{eqnarray}
where the last inequality follows from the R-RIP of $\mathcal{A}$.

The first term in (\ref{eq:lemma_estimation_proof_eq1}) has the following upper bound
\begin{eqnarray}
\norm{\mathcal{A} \left( X_0 - \mathcal{L}_{\widetilde{\Psi}} \left[ \mathcal{A} \mathcal{L}_{\widetilde{\Psi}} \right]^\dagger \mathcal{A} X_0 \right)}_2
&=& \norm{\mathcal{A} X_0 - \mathcal{A} \mathcal{L}_{\widetilde{\Psi}} \left[ \mathcal{A} \mathcal{L}_{\widetilde{\Psi}} \right]^\dagger \mathcal{A} X_0 }_2 %\nonumber \\
= \norm{\mathcal{P}_{\mathcal{R}(\mathcal{A} \mathcal{L}_{\widetilde{\Psi}})}^\perp \mathcal{A} X_0}_2 \nonumber \\
&=& \norm{\mathcal{P}_{\mathcal{R}(\mathcal{A} \mathcal{L}_{\widetilde{\Psi}})}^\perp \mathcal{A} (\mathcal{P}_{\widetilde{\Psi}} + \mathcal{P}_{\widetilde{\Psi}}^\perp) X_0}_2 %\nonumber \\
= \norm{\mathcal{P}_{\mathcal{R}(\mathcal{A} \mathcal{L}_{\widetilde{\Psi}})}^\perp \mathcal{A} (\mathcal{L}_{\widetilde{\Psi}} \mathcal{L}_{\widetilde{\Psi}}^* + \mathcal{P}_{\widetilde{\Psi}}^\perp) X_0}_2 \nonumber \\
&=& \norm{0 + \mathcal{P}_{\mathcal{R}(\mathcal{A} \mathcal{L}_{\widetilde{\Psi}})}^\perp \mathcal{A} \mathcal{P}_{\widetilde{\Psi}}^\perp X_0}_2 %\nonumber \\
\leq \|\mathcal{A} \underbrace{\mathcal{P}_{\widetilde{\Psi}}^\perp X_0}_{\rank \leq 4r}\|_2 \nonumber \\
&\leq& \sqrt{1 + \delta_{4r}(\mathcal{A})} \norm{\mathcal{P}_{\widetilde{\Psi}}^\perp X_0}_F,
\label{eq:lemma_estimation_proof_eq2}
\end{eqnarray}
where $\rank(\mathcal{P}_{\widetilde{\Psi}}^\perp X_0) \leq 4r$ holds by the subadditivity of the rank in the following way:
\begin{eqnarray*}
\rank(\mathcal{P}_{\widetilde{\Psi}}^\perp X_0)
= \rank(X_0 - \mathcal{P}_{\widetilde{\Psi}} X_0)
\leq \rank(X_0) + \rank(\mathcal{P}_{\widetilde{\Psi}} X_0)
\leq r + |\widetilde{\Psi}| \leq 4r.
\end{eqnarray*}

The second term in (\ref{eq:lemma_estimation_proof_eq1}) is bounded by
\begin{equation}
\norm{\mathcal{A} \mathcal{L}_{\widetilde{\Psi}} \left[ \mathcal{A} \mathcal{L}_{\widetilde{\Psi}} \right]^\dagger \nu}_2
= \norm{\mathcal{P}_{\mathcal{R}(\mathcal{A} \mathcal{L}_{\widetilde{\Psi}})} \nu}_2
\leq \norm{\nu}_2.
\label{eq:lemma_estimation_proof_eq3}
\end{equation}

Finally, combining (\ref{eq:lemma_estimation_proof_eq1}), (\ref{eq:lemma_estimation_proof_eq2}), and (\ref{eq:lemma_estimation_proof_eq3}) yields
\begin{eqnarray*}
\|X_0 - \widetilde{X}\|_F \leq \sqrt{\frac{1 + \delta_{4r}(\mathcal{A})}{1 - \delta_{4r}(\mathcal{A})}} \norm{\mathcal{P}_{\widetilde{\Psi}}^\perp X_0}_F + \frac{1}{\sqrt{1 - \delta_{3r}(\mathcal{A})}} \norm{\nu}_2.
\end{eqnarray*}
Applying $\delta_{3r}(\mathcal{A}) \leq \delta_{4r}(\mathcal{A}) \leq 0.04$ completes the proof.

\subsection{Proof of Lemma~\ref{lemma:prune}}
\begin{eqnarray*}
\|X_0 - \widetilde{X}_r\|_F \leq \|X_0 - \widetilde{X}\|_F + \|\widetilde{X} - \widetilde{X}_r\|_F \leq 2 \|X_0 - \widetilde{X}\|_F,
\end{eqnarray*}
where the second inequality holds by the definition of the best rank-$r$ approximation.

\subsection{Proof of Lemma~\ref{lemma:two_possibilities}}
In the $k$-th iteration, the generalized correlation maximization rule chooses $\Psi_k'$.
Let $\widetilde{\Psi}_k = \widehat{\Psi}_{k-1} \cup \Psi_k'$ in the $k$-th iteration.
Since $\widehat{\Psi}_k$ is chosen as a subset of $\widetilde{\Psi}_k$,
\begin{equation*}
\norm{ \mathcal{P}_{\widetilde{\Psi}_k}^\perp X_0}_F \leq \norm{ \mathcal{P}_{\widehat{\Psi}_k}^\perp X_0}_F, \quad \forall k \in \mathbb{Z}_+.
\end{equation*}
Then Lemma~\ref{lemma:estimation} and Lemma~\ref{lemma:prune} imply
\begin{eqnarray}
\|X_0 - \widehat{X}_k\|_F
&\leq& 2 \cdot \left(1.04 \norm{ \mathcal{P}_{\widetilde{\Psi}_k}^\perp X_0 }_F + 1.02 \norm{\nu}_2\right) \nonumber \\
&\leq& 2.08 \norm{ \mathcal{P}_{\widehat{\Psi}_k}^\perp X_0 }_F + 2.04 \norm{\nu}_2, \label{eq:ineq1_bnd_iter}
\end{eqnarray}
for all $k \in \mathbb{Z}_+$.

If $\norm{ \mathcal{P}_{\widehat{\Psi}_k}^\perp X_0}_F < 30 \norm{\nu}_2$ for some $k$, then (\ref{eq:ineq1_bnd_iter}) implies (\ref{eq:first_alt}).
The first possibility has been shown.

Otherwise, we need to show that (\ref{eq:second_alt_a}) and (\ref{eq:second_alt_b}) hold.

Assume that $\norm{ \mathcal{P}_{\widehat{\Psi}_k}^\perp X_0}_F \geq 30 \norm{\nu}_2$ for some $k \in \mathbb{Z}_+$.
Then (\ref{eq:ineq1_bnd_iter}) implies (\ref{eq:second_alt_a}).
Furthermore, we have
\begin{eqnarray*}
\|X_0 - \widehat{X}_k\|_F \geq \norm{ \mathcal{P}_{\widehat{\Psi}_k}^\perp (X_0 - \widehat{X}_k) }_F = \norm{ \mathcal{P}_{\widehat{\Psi}_k}^\perp X_0 }_F \geq 30 \norm{\nu}_2.
\end{eqnarray*}
Therefore Theorem~\ref{thm:pg_exact_low_rank} ensures that (\ref{eq:second_alt_b}) holds.

\subsection{Proof of Lemma~\ref{lemma:atom_recovery}}
Define $J$ be the set of indices for nonempty atomic bands
\begin{equation*}
J \triangleq \{ j \in \mathbb{Z}_+ :~ B_j \neq \emptyset \}.
\end{equation*}

\noindent \textbf{Claim 1} Fix an index $j \in J$.
If
\begin{equation}
\|X_0 - \widehat{X}_k\|_F \leq 2^{-(j+1)/2} \norm{X_0}_F
\label{eq:eq1_claim1}
\end{equation}
for some $k$, then
\begin{equation}
B_j \subset \widehat{\Psi}_i = \atoms{\widehat{X}_i}, \quad \forall i \geq k.
\label{eq:eq2_claim1}
\end{equation}

\begin{proof} \textbf{(Claim 1)}
First, we show that $B_j \subset \widehat{\Psi}_k$.
Assume that $B_j \not\subset \widehat{\Psi}_k$.
Then there exists $\psi \in B_j$ such that $\psi \not\in \widehat{\Psi}_k$.
This implies
\begin{equation}
\|X_0 - \widehat{X}_k\|_F \geq \norm{ \mathcal{P}_\psi (X_0 - \widehat{X}_k) }_F = \norm{ \mathcal{P}_\psi X_0 }_F > 2^{-(j+1)/2} \norm{X_0}_F,
\label{eq:cond_band_identification}
\end{equation}
which is a contradiction.
From the assumption, (\ref{eq:second_alt_b}) ensures that
\begin{equation}
\|X_0 - \widehat{X}_\ell\|_F \leq \|X_0 - \widehat{X}_k\|_F \leq 2^{-(j+1)/2} \norm{X_0}_F, \quad \forall \ell \geq k.
\label{eq:band_identification}
\end{equation}
Then Claim~1 follows.
%\qed
\end{proof}
Equation (\ref{eq:eq2_claim1}) implies that atomic band $B_j$ has been already identified in the $k$-th iteration.
From the condition of the atomic band identification in (\ref{eq:eq1_claim1}), it follows that
the identification of $B_j$ implies the identification of $B_\ell$ for all $\ell \leq j$.

\noindent \textbf{Claim 2}
Assume that $B_\ell$ has been already identified for all $\ell < j$.
Let $\beta = \left(\frac{4}{3}\right)$.
After at most
\begin{equation}
\ell = \log_\beta \left\lceil \frac{2.15 \norm{Y_j}_F}{2^{-(j+1)/2}\norm{X}_F} \right\rceil,
\label{eq:ell_claim2}
\end{equation}
more iterations, $B_j$ is identified.

\begin{proof} \textbf{(Claim 2)}
We start with the $k$-th iteration.
Since $B_\ell$ has been already identified for all $\ell < j$, $B_\ell$ is a subset of $\widehat{\Psi}_k$ for all $\ell < j$.
In other words, $\bigoplus_{\ell < j} \mathrm{span}(B_\ell) \subset \mathrm{span}(\widehat{\Psi}_k)$ and hence
\begin{equation}
\bigoplus_{\ell \geq j} \mathrm{span}(B_\ell) = \mathrm{span}(\Phi) \cap \left( \bigoplus_{\ell < j} \mathrm{span}(B_\ell) \right)^\perp \supset \mathrm{span}(\Phi) \cap \mathrm{span}(\widehat{\Psi}_k)^\perp,
\label{eq:proof_claim2a}
\end{equation}
where $\Phi = \atoms{X_0}$.
Since $X_0 \in \mathrm{span}(\Phi)$, $\mathcal{P}_{\widehat{\Psi}_k}^\perp X_0$ is the projection of $X_0$ onto $\mathrm{span}(\Phi) \cap \mathrm{span}(\widehat{\Psi}_k)^\perp$
and therefore
\begin{equation}
\norm{ Y_j }_F = \| \sum_{\ell \geq j} \mathcal{P}_{B_\ell} X_0 \|_F \geq \norm{ \mathcal{P}_{\widehat{\Psi}_k}^\perp X_0 }_F,
\label{eq:proof_claim2}
\end{equation}
where the inequality follows from (\ref{eq:proof_claim2a}).
Also note that by assumption (\ref{eq:second_alt_a}) holds.
Combining (\ref{eq:second_alt_a}) and (\ref{eq:proof_claim2}), it follows that
\begin{equation*}
\|X_0 - \widehat{X}_k\|_F \leq 2.15 \norm{ \mathcal{P}_{\widehat{\Psi}_k}^\perp X_0 }_F \leq 2.15 \norm{ Y_j }_F.
\end{equation*}
Now, $B_j$ is identified in the $(k+\ell)$-th iteration if
\begin{equation}
\|X_0 - \widehat{X}_{k+\ell}\|_F \leq 2^{-(j+1)/2} \norm{X_0}_F.
\label{eq:cond_claim2}
\end{equation}
It is easily verified that $\ell$ given in (\ref{eq:ell_claim2}) satisfies (\ref{eq:cond_claim2}).
%\qed
\end{proof}

The total number of iterations required to identify $B_j$ for all $j \in J$ is at most
\begin{equation*}
k_\star = \sum_{j \in J} \log_\beta \left\lceil 2.15 \cdot \frac{2^{(j+1)/2} \norm{Y_j}_F}{\norm{X_0}_F} \right\rceil.
\end{equation*}
For each $k \geq \lfloor k_\star \rfloor$, we have $\atoms{\widehat{X}_k} = \atoms{X_0}$.
It remains to bound $k_\star$ in terms of the profile $t = \mathrm{profile}(X_0)$.
First, note that $t = |J|$.
Using Jensen's inequality we have
\begin{eqnarray}
\exp\left\{ \frac{1}{t} \sum_{j \in J} \ln \left\lceil 2.15 \cdot \frac{2^{(j+1)/2} \norm{Y_j}_F}{\norm{X_0}_F} \right\rceil \right\}
&\leq& \exp\left\{ \frac{1}{t} \sum_{j \in J} \ln \left( 1 + 2.15 \cdot \frac{2^{(j+1)/2} \norm{Y_j}_F}{\norm{X_0}_F} \right) \right\} \nonumber \\
&\leq& \frac{1}{t} \sum_{j \in J} \left( 1 + 2.15 \cdot \frac{2^{(j+1)/2} \norm{Y_j}_F}{\norm{X_0}_F} \right)
= 1 + \frac{2.15}{t} \sum_{j \in J} \left( \frac{2^{j+1} \norm{Y_j}_F^2}{\norm{X_0}_F^2} \right)^{1/2}. \nonumber\\
\label{eq:gamean_ineq1}
\end{eqnarray}
Recall the bound on $\norm{Y_j}_F$ in (\ref{eq:bnd_norm_Yj}).
We use Jensen's inequality again and simplify the result.
\begin{eqnarray}
\frac{1}{t} \sum_{j \in J} \left( \frac{2^{j+1} \norm{Y_j}_F^2}{\norm{X_0}_F^2} \right)^{1/2}
&\leq& \frac{1}{t} \sum_{j \in J} \left( 2^{j+1} \sum_{\ell \geq j} 2^{-\ell} |B_\ell| \right)^{1/2}
\leq \left( \frac{1}{t} \sum_{j \in J} 2^{j+1} \sum_{\ell \geq j} 2^{-\ell} |B_\ell| \right)^{1/2} \nonumber \\
&\leq& \left( \frac{1}{t} \sum_{\ell \geq 0} |B_\ell| \sum_{j \leq \ell} 2^{j-\ell+1} \right)^{1/2}
\leq \left( \frac{4}{t} \sum_{\ell \geq 0} |B_\ell| \right)^{1/2}
\leq 2 \sqrt{r/t}.
\label{eq:gamean_ineq2}
\end{eqnarray}
Combining (\ref{eq:gamean_ineq1}) and (\ref{eq:gamean_ineq2}), we have
\begin{equation*}
\exp\left\{ \frac{1}{t} \sum_{j \in J} \ln \left\lceil 2.15 \cdot \frac{2^{(j+1)/2} \norm{Y_j}_F}{\norm{X_0}_F} \right\rceil \right\}
\leq 1 + 4.3 \sqrt{r/t}.
\end{equation*}
Taking the logarithm, multiplying by $t$, and dividing both sides by $\ln\beta$, we have
\begin{equation*}
k_\star \leq t \log_\beta (1 + 4.3 \sqrt{r/t}).
\end{equation*}

\subsection{Proof of Theorem~\ref{thm:bnd_iter_exact}}
Let $K = \lceil t \log_\beta (1 + 4.3 \sqrt{r/t}) \rceil$.
Suppose that (\ref{eq:first_alt}) never holds during the first $K$ iterations.
Lemma~\ref{lemma:two_possibilities} then implies that both (\ref{eq:second_alt_a}) and (\ref{eq:second_alt_b}) hold for the first $K$ iterations.
By Lemma~\ref{lemma:atom_recovery}, all atoms in $\atoms{X_0}$ are identified in the $K$-th iteration, i.e., $\widehat{\Psi}_k = \atoms{X_0}$.
Since $\widehat{\Psi}_k \subset \widetilde{\Psi}_k$, $\atoms{X_0}$ is a subset of $\widetilde{\Psi}_k$ and hence $\mathcal{P}_{\widetilde{\Psi}_k}^\perp X_0 = 0$.
Lemma~\ref{lemma:estimation} and Lemma~\ref{lemma:prune} then imply
\begin{equation*}
\|X_0 - \widehat{X}_K\|_F \leq 2 \cdot \left( 1.04 \norm{\mathcal{P}_{\widetilde{\Psi}_k}^\perp X_0}_F + 1.02 \norm{\nu}_2 \right) = 2.04 \norm{\nu}_2.
\end{equation*}
This contradicts the assumption that (\ref{eq:first_alt}) never holds during the first $K$ iterations.
Therefore, there exists $k \leq K$ where (\ref{eq:first_alt}) holds.
Repeated application of Theorem~\ref{thm:pg_exact_low_rank} gives
\begin{equation*}
\|X_0 - \widehat{X}_{K+6}\|_F < 15 \norm{\nu}_2.
\end{equation*}

\subsection{Proof of Theorem~\ref{thm:bnd_iter_gen}}
Let $X_0 \in \mathbb{C}^{m \times n}$ be an arbitrary matrix and let $t = \mathrm{profile}(X_{0,r})$.
By Lemma~\ref{lemma:noise}, we can rewrite the measurement as $b = \mathcal{A} X_{0,r} + \widetilde{\nu}$.
Theorem~\ref{thm:bnd_iter_exact} states that after at most
\begin{equation*}
t \log_{4/3}(1 + 4.3 \sqrt{r/t}) + 6
\end{equation*}
iteration, the approximation error satisfies
\begin{equation*}
\|X_{0,r} - \widehat{X}\|_F \leq 15 \norm{\widetilde{\nu}}_2.
\end{equation*}
Hence
\begin{eqnarray*}
\|X_0 - \widehat{X}\|_F
&\leq& \|X_{0,r} - \widehat{X}\|_F + \norm{X_0 - X_{0,r}}_F \\
&\leq& 15 \norm{\widetilde{\nu}}_2 + \norm{X_0 - X_{0,r}}_F \\
&\leq& 16.3 \norm{X_0 - X_{0,r}}_F + \frac{15.3}{\sqrt{r}} \norm{X_0 - X_{0,r}}_* + 15 \norm{\nu}_2 \\
&<& 17 \epsilon,
\end{eqnarray*}
where the third inequality follows from Lemma~\ref{lemma:noise}.

\bibliographystyle{IEEEtran}
\bibliography{IEEEabrv,rm}

\end{document}